\documentclass[11pt]{article}

\RequirePackage[OT1]{fontenc}

\RequirePackage[aos,amsthm,amsmath,natbib]{imsarttech}
\usepackage{amsfonts,graphicx}

\begin{document}

\begin{frontmatter}

\title{Further details on inference under right censoring 
for transformation models
with a change-point based on a covariate threshold}
\runtitle{Change-point transformation models}
\author{Michael R. Kosorok\thanksref{t1}}
\and
\author{Rui Song\thanksref{t1}}
\affiliation{University of Wisconsin-Madison}
\address{Michael R. Kosorok\\ Departments of Statistics\\
and Biostatistics \& Medical Informatics\\
1300 University Avenue\\
Madison, WI  53706\\USA\\Email: kosorok@biostat.wisc.edu}
\address{Rui Song\\ Department of Statistics\\
1300 University Avenue\\
Madison, WI  53706\\USA\\Email: rsong@stat.wisc.edu}
\runauthor{M. R. Kosorok and R. Song}
\thankstext{t1}{Supported in part by Grant CA075142 from the 
National Cancer Institute.}

\begin{abstract}
We consider linear transformation models applied to right censored survival
data with a change-point in the regression coefficient based on a covariate
threshold.  We establish consistency and weak convergence of the
nonparametric maximum likelihood estimators.  The change-point parameter
is shown to be $\,n$-consistent, while the remaining parameters are shown to 
have the expected root-$n$ consistency. We show that the 
procedure is adaptive in the sense that the non-threshold parameters are 
estimable with the same precision as if the true threshold value were known.
We also develop Monte-Carlo methods of inference for model parameters
and score tests for the existence of a change-point.
A key difficulty here is that some of the model
parameters are not identifiable under the null hypothesis of no change-point.
Simulation studies establish the validity of the proposed 
score tests for finite sample sizes.
\end{abstract}



\end{frontmatter}

\newtheorem{theorem}{\indent \sc Theorem}
\newtheorem{corollary}{\indent \sc Corollary}
\newtheorem{lemma}{\indent \sc Lemma}
\newtheorem{proposition}{\indent \sc Proposition}
\newtheorem{remark}{\indent \sc Remark}
\newcommand{\phif}{\textsc{igf}}
\newcommand{\sign}{\mbox{sign}}
\newcommand{\phgf}{\textsc{gf}}
\newcommand{\fix}{$\textsc{gf}_0$}
\newcommand{\mb}[1]{\mbox{\bf #1}}
\newcommand{\Exp}[1]{\mbox{E}\left[#1\right]}
\newcommand{\pr}[1]{\mbox{P}\left[#1\right]}
\newcommand{\pp}[0]{\mathbb{P}}
\newcommand{\ee}[0]{\mbox{E}}
\newcommand{\re}[0]{\mathbb{R}}
\newcommand{\argmax}[0]{\mbox{argmax}}
\newcommand{\argmin}[0]{\mbox{argmin}}
\newcommand{\ind}[0]{\mbox{\Large\bf 1}}
\newcommand{\narrow}{\stackrel{n\rightarrow\infty}{\longrightarrow}}
\newcommand{\weakpn}{\stackrel{P_n}{\leadsto}}
\newcommand{\weakpnboot}{\mbox{\raisebox{-1.5ex}{$\stackrel
{\mbox{\scriptsize $P_n$}}{\stackrel{\mbox{\normalsize $\leadsto$}}
{\mbox{\normalsize $\circ$}}}$}}\,}
\newcommand{\ol}[1]{\overline{#1}}
\newcommand{\avgse}[1] { \bar{\hat{\sigma}}_{#1} }
\newcommand{\mcse}[1]  { \sigma^{*}_{#1} }
\newcommand{\po}{\textsc{po}}
\newcommand{\ph}{\textsc{ph}}

\section{Introduction} The linear transformation model states that
a continuous outcome $U$, given a $d$-dimensional covariate vector $Z$,
has the form
\begin{eqnarray}
 H(U)= - \beta ^\prime Z + \varepsilon, \label{s1.e1}
\end{eqnarray}
where $H$ is an increasing, unknown transformation function, $\beta\in\re^d$ 
are the unknown regression parameters of interest, and $\varepsilon$ has a 
known distribution $F$. This model is readily applied to a failure time $T$
by letting $U=\log T$ and $H(u)=\log A(e^u)$, where $A$ is an
unspecified integrated baseline hazard. Setting 
$F(s)=1-\exp(-e^s)$ results in the Cox model, while setting 
$F(s)=e^s/(1+e^s)$ results in the proportional odds model. More
generally, the transformation model for a survival time $T$ conditionally
on a time-dependent covariate $\tilde{Z}(t)=\{Z(s),0\leq s\leq t\}$, 
takes the form 
\begin{eqnarray}
\pr{T>t|\tilde{Z}(t)}=S_Z(t)&\equiv&\Lambda\left(\int_0^te^{\beta'Z(s)}
dA(s)\right),\label{new.e1}
\end{eqnarray}
where $\Lambda$ is a known decreasing function
with $\Lambda(0)=1$. The model~(\ref{new.e1})
becomes model~(\ref{s1.e1}) when the covariates are time-independent
and $F(s)=1-\Lambda(e^s)$.

In data analysis, the assumption of linearity of the regression
effect in~(\ref{new.e1}) is not always satisfied over the whole
range of the covariate, and the fit may be improved with a
two-phase transformation model having a change-point at an
unknown threshold of a one-dimensional covariate $Y$. Let
$Z=(Z_1,Z_2)$, where $Z_1$ and $Z_2$ are possibly time-dependent
covariates in $\re^p$ and $\re^q$, respectively, where $p+q=d$
and $q\geq 1$. The new model is obtained by replacing $\beta'Z(s)$ 
in~(\ref{new.e1}) with
\begin{eqnarray}
r_{\xi}(s;Z,Y)\equiv\beta'Z(s) +  [\alpha + \eta'Z_2(s)]\ind\{Y >
\zeta\},\label{new.e2}
\end{eqnarray}
where $\alpha$ is a scalar, $\eta\in\re^q$, $\ind\{B\}$ is the
indicator of $B$, and $\xi$ 
denotes the collected parameters $(\alpha,\beta,\eta,\zeta)$. 
We also require $Y$ to be time-independent but allow it to possibly
be one of the covariates in $Z(t)$. The overall goal of this paper
is to develop methods of inference for this model applied to
right censored data.

We note that for the special case when $\alpha=0$ and 
$\Lambda(t)=e^{-t}$, the model~(\ref{new.e2})
becomes the Cox model considered by \cite{p03} under a slightly 
different parameterization. Permitting a nonzero $\alpha$
allows the possibility of a ``bent-line''
covariate effect. Suppose, for example, that 
$Z_2$ is one-dimensional and time-independent, 
while $Z_1\in\re^{d-1}$ may be time-dependent.
If we set $Y=Z_2$ and $\beta=(\beta_1',\beta_2')'$, where
$\beta_1\in\re^{d-1}$ and $\beta_2\in\re$, the model~(\ref{new.e2})
becomes $r_{\xi}(s;Z,Y)=\beta_1'Z_1(s)+\beta_2Z_2
+(\alpha+\eta Z_2)\ind\{Z_2>\zeta\}$. 
When $\alpha=-\eta\zeta$, the covariate effect for $Z_2$
consists of two connected linear segments.  In many biological settings, 
such a bent-line effect is realistic and can be much easier
to interpret than a quadratic or more complex nonlinear effect \cite{c89}. 
Hence including the intercept term $\alpha$ is useful
for applications.

Linear transformation models of the form~(\ref{s1.e1}) have been
widely used and studied (see, for example, \cite{bc64,bd81,bc82,p82,dd88,
cwy95,cwy97,fyw98,bn04}). Efficient methods of estimation in the 
uncensored setting were rigorously studied by \cite{br97}, 
among others. The model~(\ref{new.e1}) for right-censored data has 
also been studied rigorously for a variety of specific choices of 
$\Lambda$ \cite{p84,mrv97,stg98,s98}; for
general but known $\Lambda$ \cite{sv04}; and for certain 
parameterized families of $\Lambda$ \cite{klf04}.

Change-point models have also been studied extensively and have
proven to be popular in clinical research. Several researchers have
considered a nonregular Cox model involving a two-phase regression on
time-dependent covariates, with a change-point at an unknown time
\cite{lsl90,lb97,ltc97}. As mentioned above, \cite{p03}
considered the Cox model with a
change-point at an unknown threshold of a covariate.
These authors studied the maximum partial likelihood estimators of
the parameters and the estimator of the baseline hazard function.
They show that the estimator of the threshold parameter
is $n$-consistent, while the regression parameters are 
$\sqrt{n}$-consistent. This happens because
the likelihood function is not differentiable with respect to the 
threshold parameter, and hence the usual Taylor expansion
is not available. In this paper, we focus on
the covariate threshold setting. While time threshold models are
also interesting, we will not pursue them further in this paper
because the underlying techniques for estimation and
inference are quite distinct from the covariate threshold setting. 

The contribution of our paper builds on \cite{p03} in three
important ways.  Firstly, we extend to general transformation models.
This results in a significant increase in complexity over the Cox
model since estimation of the baseline hazard can no longer be
avoided through the use of the partial-profile likelihood.  Secondly,
we study nonparametric maximum likelihood inference for all model parameters.
As part of this, we show that the estimation procedure is adaptive in the 
sense that the non-threshold parameters---including the infinite-dimensional
parameter~$A$---are estimable with the same precision 
as if the true threshold parameter were known. Thirdly, we develop hypothesis
tests for the existence of a change-point. This is quite challenging
since some of the model parameters are no longer identifiable under
the null hypothesis of no change-point.  \cite{a01} considers
similar nonstandard testing problems when the model is fully
parametric and establishes asymptotic null and local alternative 
distributions of a number of likelihood-based test procedures. 
Unfortunately, Andrews' results are not directly
applicable to our setting because of the presence of an infinite 
dimensional nuisance parameter, the baseline integrated hazard $A$, 
and new methods are required.

The next section, section~2, presents the data and model
assumptions. The nonparametric maximum log-likelihood estimation 
(NPMLE) procedure is
presented in section~3. In section~4, we establish the consistency
of the estimators. Score and information operators of the regular
parameters are given in section~5. Results on the convergence rates
of the estimators are established in section~6. Section~7 presents
weak convergence results for the estimators, including the
asymptotic distribution of the change-point estimator and
the asymptotic normality of the other parameters. This section
also establishes the adaptive semiparametric efficiency
mentioned above. Monte Carlo inference for the parameters is discussed
in section~8. Methods for testing the existence of a change-point are
then presented in section~9. A brief discussion on implementation
and a small simulation study evaluating the moderate
sample size performance of the proposed change-point tests are
given in section~10. Proofs are given in section~11.

\section{The data set-up and model assumptions}
The data $X_i=(V_i,\delta_i,$ $Z_i,Y_i)$, $i=1,\ldots,n$, consists of $n$ 
i.i.d. realizations of $X = (V,\delta,Z,Y)$, where $V = T \land C$, $\delta
= 1(T \le C)$, and $C$ is a right censoring time. The analysis is
restricted to the interval $[0, \tau]$, where $\tau < \infty$.
The covariate $Y\in\re$ and $Z \equiv \{Z(t), t \in [0, \tau] \}$ is assumed 
to be a caglad (left-continuous with right-hand limits) process with
$Z(t)=(Z_1'(t),Z_2'(t))'\in\re^p\times\re^q$, for all $t\in [0, \tau]$, 
where $q\geq 1$ but $p=0$ is allowed.

We assume that conditionally on $Z$ and $Y$, the survival function
at time $t$ has the form:
\begin{eqnarray}
S_{Z,Y}(t)\equiv\Lambda \left(\int_0^t e^{r_{\xi}(u;Z,Y)}dA(u)\right),
\label{s2.e1}
\end{eqnarray}
where $\Lambda$ is a known, thrice differentiable 
decreasing function with $\Lambda(0)=1$,
$r_{\xi}(s;Z,Y)$ is as defined in~(\ref{new.e2}), and $A$ is an
unknown increasing function restricted to $[0,\tau]$. 

Let $G\equiv-\log\Lambda$, 
and define the derivatives 
$\dot{\Lambda}\equiv\partial\Lambda(t)/(\partial t)$,
$\ddot{\Lambda}\equiv\partial\dot{\Lambda}(t)/(\partial t)$,
$\dot G \equiv \partial G(t) /(\partial t)$,
$\ddot G \equiv \partial \dot G (t) /(\partial t)$, and 
$\dddot G \equiv \partial\ddot G /(\partial t)$.
We also define the collected parameters
$\gamma\equiv(\alpha,\eta,\beta)$, $\psi\equiv(\gamma, A)$, and
$\theta \equiv(\psi, \zeta)$. We use $P$ to denote the
true probability measure, while the true parameter values are
indicated with a subscript 0.

We now make the following additional assumptions:
\begin{itemize}
\item[A1]: $P[C=0]=0$, $P[C \ge \tau | Z,Y] = P[C = \tau | Z,Y] > 0$
almost surely, and censoring is independent of $T$ given $(Z,Y)$
and uninformative. 
\item[A2]: The total variation of $Z(\cdot)$ on $[0, \tau]$ is 
$\le m_0<\infty$ almost surely.
\item[B1]: $\zeta_0\in(a,b)$, for some known $-\infty<a<b<\infty$
with $P[Y<a]>0$ and $P[Y>b]>0$.
\item[B2]: For some neighborhood $\tilde{V}(\zeta_0)$ of $\zeta_0$:
\begin{itemize}
\item[(i)] the density of $Y$, $\tilde{h}$, exists and
is strictly positive, bounded and continuous for all 
$y\in\tilde{V}(\zeta_0)$; and
\item[(ii)] the conditional law of $(C,Z)$ given $Y=y$, 
${\cal L}_y$, is left-continuous with right-hand limits 
over $\tilde{V}(\zeta_0)$.
\end{itemize}
\item[B3]: For some $t_1,t_2\in(0,\tau]$, both var$[Z(t_1)|Y=\zeta_0]$ 
and var$[Z(t_2)|Y=\zeta_0+]$ are positive definite.
\item[B4]: For some $t_3,t_4\in(0,\tau]$, both
var$[Z(t_3)|Y<a]$ and var$[Z(t_4)|Y>b]$ are positive definite.
\item[C1]: $\alpha_0\in\Upsilon\subset\re$, $\beta_0\in B_1\subset\re^d$,
$\eta_0\in B_2\subset\re^q$, where $d\geq q\geq 1$, and $\Upsilon$, $B_1$
and $B_2$ are open, convex, bounded and known.
\item[C2]: Either $\alpha_0\neq 0$ or $\eta_0\neq 0$.
\item[C3]: $A_0\in{\cal A}$, where ${\cal A}$ is the set of all
increasing functions $A:[0,\tau]\mapsto[0,\infty)$ with
$A(0)=0$ and $A(\tau)<\infty$; and $A_0$ has derivative $a_0$ 
satisfying $0<a_0(t)<\infty$ for all $t\in[0,\tau]$.
\item[D1]: $G:[0,\infty)\mapsto[0,\infty)$ is thrice continuously
differentiable, with $G(0)=0$, and, for each $u\in[0,\infty)$, 
$0<\dot{G}(u),\ddot{\Lambda}(u)<\infty$ and 
$\sup_{s\in[0,u]}|\dddot G(s)|<\infty$. 
\item[D2]: For some $c_0>0$, both 
$\sup_{u\geq 0}|u^{c_0}\Lambda(u)|<\infty$ and
$\sup_{u\geq 0}|u^{1+c_0}\dot{\Lambda}(u)|<\infty$.
\end{itemize}

Conditions~A1, A2, C1 and~C3 are commonly used for NPMLE
consistency and identifiability in right-censored
transformation models, while conditions~B1, B2, B3 and~C2 are
needed for change-point identifiability. As pointed out by a
referee, the use of a time-dependent covariate will require
that $Z_i(V_j)$ be observed for each individual $i$ and for
every $j$ such that $\delta_1=1$ and $V_j\leq V_i$. While this
is often assumed in theoretical contexts, it can be unrealistic
in practice, where missing values of $Z_i(t)$ are not
unusual (see \cite{ly93}).  Frequently, data analysts will simply
carry the last observation of $Z_i(t)$ forward to avoid the missingness
problem. Unfortunately, this simple solution is not necessarily valid.
However, addressing this
issue thoroughly is beyond the scope of this paper, and we will
only mention it again briefly in section~9, 
where we develop a test of the null hypothesis that there is no
change-point ($H_0:\alpha_0=0$ and $\eta_0=0$). Also in section~9,
we will relax condition~C2 
to allow for a sequence of contiguous alternative hypotheses that 
includes $H_0$. Condition~B2(ii) is also needed to obtain weak convergence
for the NPMLE of $\zeta_0$. The continuity requirements 
at each point $y$ can be restated in the following way: 
${\cal L}_{\zeta}$ converges
weakly to ${\cal L}_y$, as $\zeta\uparrow y$; and 
${\cal L}_{\zeta}$ converges weakly to ${\cal L}_{y+}$, 
as $\zeta\downarrow y$, for some law ${\cal L}_{y+}$.
It would require a fairly pathological relationship among
the variables $(C,Z,Y)$ for this not to hold. Condition~B4 will also be 
needed for the change-point test developed in section~9. 

Conditions~D1 and~D2 
are also needed for asymptotic normality. Condition D1~is quite similar to
conditions~(G.1) through~(G.4) in \cite{sv04} who
use the condition for developing asymptotic theory
for transformation models without a change-point. Condition~D2
is slightly weaker than conditions~D2 and~D3 of \cite{klf04}
who use the condition to obtain asymptotic theory
for frailty regression models without a change-point.
The following are several instances that satisfy conditions~D1 and D2:
\begin{enumerate}
\item $\Lambda(u)=e^{-u}$ corresponds to the extreme value distribution
and results in the Cox model. 
\item $\Lambda(u)=(1+c u)^{-1/c}$, for any $c\in(0,\infty)$,
corresponds to the family of log-Pareto distributions and results
in the odds-rate transformation family. Taking the limit as 
$c\downarrow 0$ yields the Cox model, while $c=1$ 
yields the proportional odds model. 
\item $\Lambda(u)=\Exp{e^{-Wu}}$, where $W$ is a positive frailty with
$\Exp{W^{-c}}<\infty$, for some $c>0$, and 
$\Exp{W^4}<\infty$, corresponds to the family of frailty transformations.
In addition to the odds-rate family, these conditions are
satisfied by both the inverse Gaussian and log-normal families
(see \cite{klf04}), as well as many other frailty families.
\item $\Lambda(u)=[1+2cu+u^2]^{-1}$, 
where $c\in(1/2,1)$. Because this is the Laplace transform of
$t\mapsto e^{-ct}$ $\times\sin\left(t\sqrt{1-c^2}\right)/\sqrt{1-c^2}$, 
it is not the Laplace transform of a density. Hence this family 
is not a member of the family of frailty transformations. Note, however, that
taking the limit as $c\uparrow 1$ results in the Laplace transform
of the frailty density $te^{-t}$.
\end{enumerate}

Verification of these conditions is routine for examples~1, 2 and~4 above, but
verification for example~3 is slightly more involved:
\begin{lemma}\label{l.v1}
Conditions~D1 and~D2 are satisfied for example~3 above.
\end{lemma}

\section{Nonparametric Maximum log-likelihood estimation} The
nonparametric log-likelihood has the form $L_n(\psi,\zeta)\equiv$
\begin{eqnarray}
&&\mathbb{P}_n \left\{ 
\delta\log(a(V))+l_1^{\psi}(V,\delta,Z) 
\ind\{Y \le \zeta \} + l_2^{\psi}(V,\delta,Z) 
\ind\{Y > \zeta \} \right \}, \label{s3.e1}
\end{eqnarray}
where
\begin{eqnarray*}
l_1^{\psi}(V,\delta,Z)&\equiv&\int_0^{\tau}\left[\log\dot{G}
\left(H^{\psi}_1(s)\right)+\beta' Z(s)\right]
dN(s)-G(H^{\psi}_1(V)),\\
l_2^{\psi}(V,\delta,Z)
&\equiv&\int_0^{\tau}\left[\log\dot{G}\left(H^{\psi}_2(s)
\right)+\beta'Z(s)+\alpha+\eta'Z_2(s)\right]dN(s)\\
&&-G(H^{\psi}_2(V)),
\end{eqnarray*}
where $N(t)\equiv\ind\{V\leq t\}\delta$, $\tilde{Y}(s)\equiv\ind\{V\geq s\}$,
$a \equiv dA/dt$, 
$H^{\psi}_1(t)\equiv\int_0^t\tilde{Y}(s)e^{\beta' Z(s)}dA(s)$,
$H^{\psi}_2(t)\equiv\int_0^t\tilde{Y}(s)e^{\beta'Z(s)+\alpha + \eta'Z_2(s)}
dA(s)$, and $\mathbb{P}_n$ is the
empirical probability measure.

As discussed by \cite{mrv97}, the
maximum likelihood estimator for $a$ does not exist, since any
unrestricted maximizer of~(\ref{s3.e1}) puts mass only at observed failure
times and is thus not a continuous hazard. We replace $a(u)$ 
in $L_n(\psi,\zeta)$ with $n\Delta A(u)$ as suggested in \cite{p98}
who remarked that
this form of the empirical log-likelihood function is asymptotically
equal to the true log-likelihood function in certain instances. 
Let $\tilde{L}_n(\psi,\zeta)$ be this modified log-likelihood.
Note that the maximum likelihood
estimator for $\zeta$ is not unique, since the likelihood is constant
in $\zeta$ over the intervals $[Y_{(r)},Y_{(r+1)})$, where $Y_{(1)}  
<\cdots<Y_{(r)}<\cdots<Y_{(n)}$ are the order statistics of~$Y$. For this
reason, we only need to consider $\zeta$ at the values of the
$Y$ order statistics.

The estimators are obtained in the following way: For fixed
$\zeta$, we maximize the fully nonparametric log-likelihood over
$\psi$, to obtain the profile log-likelihood
$pL_n(\zeta)\equiv\sup_{\psi}\tilde{L}_n(\psi,\zeta)$. We then maximize
$pL_n(\zeta)$ over $\zeta$, to obtain $\hat{\zeta}_n$; and then
compute $\hat{\psi}_n=\argmax_{\psi}\tilde{L}_n(\psi,\hat{\zeta}_n)$.
This yields the NPMLE $\hat{\theta}_n=(\hat{\psi}_n,\hat{\zeta}_n)$
for $\theta_0$. Hence we obtain an estimator for $A_0$ but not for $a_0$.

\section{Consistency}
To study consistency, we first 
characterize the NPMLE $\hat\theta_n$.
Consider the following one-dimensional submodels for A:
\begin{eqnarray*}
t \mapsto A_t \equiv \int_0 ^{(\cdot)} (1 + tg(s) )dA(s),
\end{eqnarray*}
where $g$ is an arbitrary non-negative bounded function. A score
function for $A$, defined as the derivative of $\tilde{L}_n(\xi, A_t)$
with respect to $t$ at $t=0$, is
\begin{eqnarray}
&&\mathbb{P}_n \left \{ \delta g(X) - \left[ 
\dot{G}(H^{\theta}(V))-
\delta\frac{ \ddot{G}(H^{\theta}(V))}
{ \dot{G}(H^{\theta}(V))}\right] 
\int_0^{\tau}\tilde{Y}(s)
e^{r_{\xi}(s;Z,Y)}g(s)dA(s) \right\}, ~\label{c4:e1}
\end{eqnarray}
where $H^{\theta}(t)\equiv\int_0^t\tilde{Y}(s)
e^{r_{\xi}(s;Z,Y)}dA(s)$.
For any fixed $\xi$, let $\hat A_{\xi}$ denote the maximizer of
$A\mapsto\tilde{L}_n(\xi, A)$, and let $\hat \theta_{\xi} \equiv (\xi, \hat
A_{\xi})$. Then the score function~(\ref{c4:e1}) is equal to zero
when evaluated at $\hat \theta_{\xi}$. We select $g(u) = \ind
\{u \le t \}$, insert this into~(\ref{c4:e1}), and equate the
resulting expression to zero: $\hat{A}_{\xi}(u)=$
\begin{eqnarray}
\label{c4:e2}&&\\
\int_0^u\left(\pp_n \left[\tilde{Y}(s)e^{r_{\xi}(s;Z,Y)}\left(
\dot{G}\left\{H^{\hat{\theta}_{\xi}}(V)\right\}-\delta\frac{
\ddot{G}\left\{H^{\hat{\theta}_{\xi}}(V)\right\}}{\dot{G}
\left\{H^{\hat{\theta}_{\xi}}(V)\right\}}
\right)\right]\right)^{-1}\pp_n\{dN(s)\}&&\nonumber\\
\equiv\int_0^u\{\mathbb{P}_nW(s;
\hat{\theta}_{\xi})\}^{-1}\pp_n\{dN(s)\}.&&\nonumber
\end{eqnarray}
Now the profile likelihood has the form
$pL_n(\zeta)=\argmax_{\gamma}\tilde{L}_n
\left((\gamma,\hat{A}_{(\gamma,\zeta)}),\zeta\right)$.

The above characterization facilitates the following consistency
results for $\hat{\theta}_n$:
\begin{lemma}\label{l1}
Under the regularity conditions of section~2,
the transformation model with a change-point based on a covariate
threshold is identifiable.
\end{lemma}
\begin{lemma}\label{l2}
Under the regularity conditions of section~2,
$\hat{A}_n$ is asymptotically bounded, and thus the NPMLE
$\hat{\theta}_n$ exists.
\end{lemma}
Using these results, we can establish the uniform consistency of
$\hat \theta_n$:
\begin{theorem}\label{t1}
Under the regularity conditions of section~2, 
$\hat {\theta}_n$ converges outer almost surely to $\theta_0$ 
in the uniform norm.
\end{theorem}

\section{Score and information operators for regular parameters}
In this section, we derive the score and information operators
for the collected parameters $\psi$. We refer to these parameters
as the regular parameters because, as we will see in section~6,
these parameters converge at the $\sqrt{n}$ rate. On the other
hand, $\hat{\zeta}_n$ converges at the $n$ rate and thus 
the parameter $\zeta$ is not regular. The score and
information operators for $\psi$ are needed for the convergence
rate and weak limit results of sections~6 and~7. 

Let $\mathcal{H}$ denote the space of the elements $h = (h_1, h_2, h_3, h_4)$ 
such that $h_1 \in \mathbb{R}$, $h_2\in\re^q$, $h_3 \in \mathbb{R}^d$, 
and $h_4 \in D[0,\tau]$, where $D[0,\tau]$ is the space of cadlag
functions (right-continuous with left-hand limits) on $[0,\tau]$. 
We denote by $BV$ the subspace of $D[0,\tau]$ consisting of functions 
that are of bounded variation over the 
interval $[0,\tau]$. Define, for future use, 
the following linear functional for 
each $\theta=(\psi,\zeta)$ and each $t \in [0,\tau]$:
\begin{eqnarray}
R^t_{\zeta,\psi}(f) \equiv \int_0^t f(u)\tilde{Y}(u)
e^{r_{\xi}(u;Z,Y)}dA(u),
\end{eqnarray}
where $f$ is an element or vector of elements in $BV$.  
Also let $\rho_1 (h) \equiv ( | h_1 |^2 +  
\|h_2 \|^2+ \| h_3\|^2+  \| h_4 \|_v^2)^{1/2}$ and  $\mathcal{H}_r
\equiv \{h \in \mathcal{H}: \rho_1(h)\leq r \}$, where $\|\cdot\|_v$
is the total variation norm on $BV$ and $r\in(0,\infty)$. 

The parameter
$\psi\in\Psi\equiv\Upsilon\times B_2\times B_1\times{\cal A}$ 
can be considered a linear functional on $\mathcal{H}_r$ by defining
$\psi(h) \equiv h_1\alpha + h_2' \eta + h_3' \beta + \int_0^{\tau}
h_4(u)dA(u)$, $h \in \mathcal{H}_r$.
Viewed this way, $\Psi$ is a subset of $\ell^{\infty}(H_r)$ with uniform
norm $\|\psi\|_{(r)}\equiv\sup_{h\in{\cal H}_r}|\psi(h)|$, where
$\ell^{\infty}(B)$ is the space of bounded functionals on~$B$. Note that
${\cal H}_1$ is rich enough to extract all components of $\psi$. This
is easy to see for the Euclidean components; and, for the $A$ component, 
it works by using the elements $\{h:h_1=0,h_2=0,h_3=0,
h_4(u)=\ind\{u\leq t\},t\in[0,\tau]\}\subset{\cal H}_1$. 

In section~5.1, we derive the score operator; while in section~5.2 
we derive the information operator and establish its continuous invertibility.

\subsection{The score operator} Using the one-dimensional submodel
\begin{eqnarray*}
 t \rightarrow \psi_t \equiv \psi + t(h_1, h_2, h_3, 
\int_0^{(\cdot)}h_4(u)dA(u)), ~~~ h \in \mathcal{H}_r,
\end{eqnarray*}
the score operator takes the form
\[U^{\tau}_{n\zeta}(\psi)(h) \equiv
\left.\frac{\partial}{\partial t}L_n(\psi_t, \zeta)  \right|_{t=0}
=\mathbb{P}_n U^{\tau}_{\zeta}(\psi)(h),\]
where $U_{\zeta}^{\tau}(\psi)(h)\equiv U^{\tau}_{\zeta,1}(\psi)(h_1) +
U^{\tau}_{\zeta,2}(\psi)(h_2)+
U^{\tau}_{\zeta,3}(\psi)(h_3)+U^{\tau}_{\zeta,4}(\psi)(h_4)$, and
\begin{eqnarray*}
 U^{\tau}_{\zeta,1}(\psi)(h_1)&\equiv&
\ind\{Y>\zeta\}\left\{\int_0^{\tau}h_1dN(u)-\hat{\Xi}_{\theta}^{(0)}(\tau)
R^{\tau}_{\zeta,\psi}(h_1)\right\},\\
U^{\tau}_{\zeta,2}(\psi)(h_2)&\equiv&\ind(Y>\zeta)\left\{
\int_0^{\tau}Z_2'(u)h_2dN(u)-\hat{\Xi}_{\theta}^{(0)}(\tau)
R^{\tau}_{\zeta,\psi}(Z_2'h_2)\right\},\\
U^{\tau}_{\zeta,3}(\psi)(h_3)&\equiv& 
\int_0^{\tau}Z'(u)h_3dN(u)-\hat{\Xi}_{\theta}^{(0)}(\tau)
R^{\tau}_{\zeta,\psi}(Z'h_3),\\
U^{\tau}_{\zeta,4}(\psi)(h_4)&\equiv&\int_0^{\tau}h_4(u)dN(u) 
-\hat{\Xi}_{\theta}^{(0)}(\tau)R^{\tau}_{\zeta,\psi}(h_4),\\
\hat{\Xi}_{\theta}^{(0)}(\tau)&\equiv&\ind\{Y\leq\zeta\}
\hat{\Xi}_{\psi,1}^{(0)}(\tau)+\ind\{Y>\zeta\}
\hat{\Xi}_{\psi,2}^{(0)}(\tau),
\end{eqnarray*}
and where, for $j=1,2$,
\[\hat{\Xi}_{\psi,j}^{(0)}(\tau)\equiv
\left[\dot{G}(H^{\psi}_j(V \wedge \tau))  -
\delta \frac{ \ddot{G}(H^{\psi}_j(V \wedge \tau))}
{\dot{G}(H^{\psi}_j(V \wedge \tau))}\right].\]
The dependence in the notation on $\tau$ will prove useful in later
developments.

\subsection{The information operator} To obtain the information
operator, we can differentiate the expectation of the score operator
using the map $t \rightarrow \psi+ t\psi_1$,
where $\psi,\psi_1\in\Psi$. The information operator,
$\sigma_{\theta}:\mathcal{H}_{\infty} \rightarrow
\mathcal{H}_{\infty}$, where $\mathcal{H}_{\infty}\equiv\{h:\mbox{$h\in
\mathcal{H}_r$ for some $r<\infty$}\}$, satisfies
\begin{eqnarray}
\psi_1(\sigma_{\theta}(h))
&=&\left. -\frac{\partial}{\partial t}
PU^{\tau}_{\zeta}(\psi+t\psi_1)(h) \right|_{t=0},\label{new.j12.e1}
\end{eqnarray}
for every $h\in{\cal H}_{\infty}$. 
Taking the G\^{a}teaux derivative in~(\ref{new.j12.e1}), we obtain 
$\sigma_{\theta}(h)=$
\begin{eqnarray}
\label{c5.e2}&&\\
\left(\begin{array}{cccc}\sigma_{\theta}^{11}&
\sigma_{\theta}^{12}& \sigma_{\theta}^{13} & \sigma_{\theta}^{14} \\
\sigma_{\theta}^{21}&\sigma_{\theta}^{22}&
\sigma_{\theta}^{23}& \sigma_{\theta}^{24} \\
\sigma_{\theta}^{31}&\sigma_{\theta}^{32}
&\sigma_{\theta}^{33}&\sigma_{\theta}^{34}\\
\sigma_{\theta}^{41}&\sigma_{\theta}^{42}&\sigma_{\theta}^{43}
&\sigma_{\theta}^{44}
\end{array}\right)
\left(\begin{array}{c}h_1\\ h_2\\ h_3\\ h_4\end{array}\right)
\equiv P\left(\begin{array}{cccc}\hat{\sigma}_{\theta}^{11}&
\hat{\sigma}_{\theta}^{12}& \hat{\sigma}_{\theta}^{13}& 
\hat{\sigma}_{\theta}^{14} \\
\hat{\sigma}_{\theta}^{21}&
\hat{\sigma}_{\theta}^{22}&\hat{\sigma}_{\theta}^{23}& 
\hat{\sigma}_{\theta}^{24} \\
\hat{\sigma}_{\theta}^{31}&\hat{\sigma}_{\theta}^{32}&
\hat{\sigma}_{\theta}^{33}& \hat{\sigma}_{\theta}^{34}\\
\hat{\sigma}_{\theta}^{41}&\hat{\sigma}_{\theta}^{42}&
\hat{\sigma}_{\theta}^{43}& \hat{\sigma}_{\theta}^{44}
\end{array}\right)
\left(\begin{array}{c}h_1\\ h_2\\ h_3\\ h_4\end{array}\right)&&
\nonumber
\end{eqnarray}
$\equiv P\hat{\sigma}_{\theta}(h)$, where
\begin{eqnarray*}
\hat{\sigma}_{\theta}^{11}(h_1)&\equiv&\ind\{Y>\zeta\}\left\{
\hat{\Xi}_{\theta}^{(0)}(\tau)+\hat{\Xi}_{\theta}^{(1)}(\tau)
H_2^{\psi}(V\wedge\tau)\right\}R_{\zeta,\psi}^{\tau}(h_1),\\
\hat{\sigma}_{\theta}^{12}(h_2)&\equiv&\ind\{Y>\zeta\}
\left\{\hat{\Xi}_{\theta}^{(0)}(\tau)+\hat{\Xi}_{\theta}^{(1)}(\tau)
H_2^{\psi}(V\wedge\tau)\right\}R_{\zeta,\psi}^{\tau}(Z_2'h_2),\\
\hat{\sigma}_{\theta}^{13}(h_3)&\equiv&\ind\{Y>\zeta\}
\left\{\hat{\Xi}_{\theta}^{(0)}(\tau)+\hat{\Xi}_{\theta}^{(1)}(\tau)
H_2^{\psi}(V\wedge\tau)\right\}R_{\zeta,\psi}^{\tau}(Z'h_3),\\
\hat{\sigma}_{\theta}^{14}(h_4)&\equiv&\ind\{Y>\zeta\}
\left\{\hat{\Xi}_{\theta}^{(0)}(\tau)+\hat{\Xi}_{\theta}^{(1)}(\tau)
H_2^{\psi}(V\wedge\tau)\right\}R_{\zeta,\psi}^{\tau}(h_4),\\
\hat{\sigma}_{\theta}^{21}(h_1)&\equiv&\ind\{Y>\zeta\}\left\{
\hat{\Xi}_{\theta}^{(0)}(\tau)R_{\zeta,\psi}^{\tau}(Z_2h_1)
+\hat{\Xi}_{\theta}^{(1)}(\tau)R_{\zeta,\psi}^{\tau}(Z_2)
R_{\zeta,\psi}^{\tau}(h_1)\right\}\\
\hat{\sigma}_{\theta}^{22}(h_2)&\equiv&\ind\{Y>\zeta\}\left\{
\hat{\Xi}_{\theta}^{(0)}(\tau)R_{\zeta,\psi}^{\tau}(Z_2Z_2'h_2)
+\hat{\Xi}_{\theta}^{(1)}(\tau)R_{\zeta,\psi}^{\tau}(Z_2)
R_{\zeta,\psi}^{\tau}(Z_2'h_2)\right\},\\
\hat{\sigma}_{\theta}^{23}(h_3)&\equiv& 
\hat{\Xi}_{\theta}^{(0)}(\tau)R_{\zeta,\psi}^{\tau}(Z_2Z'h_3)
+\hat{\Xi}_{\theta}^{(1)}(\tau)R_{\zeta,\psi}^{\tau}(Z_2)
R_{\zeta,\psi}^{\tau}(Z'h_3),\\
\hat{\sigma}_{\theta}^{24}(h_4)&\equiv& 
\hat{\Xi}_{\theta}^{(0)}(\tau)R_{\zeta,\psi}^{\tau}(Z_2h_4)
+\hat{\Xi}_{\theta}^{(1)}(\tau)R_{\zeta,\psi}^{\tau}(Z_2)
R_{\zeta,\psi}^{\tau}(h_4),\\
\hat{\sigma}_{\theta}^{31}(h_1)&\equiv&\ind\{Y>\zeta\}\left\{
\hat{\Xi}_{\theta}^{(0)}(\tau)R_{\zeta,\psi}^{\tau}(Zh_1)
+\hat{\Xi}_{\theta}^{(1)}(\tau)R_{\zeta,\psi}^{\tau}(Z)
R_{\zeta,\psi}^{\tau}(h_1)\right\},\\
\hat{\sigma}_{\theta}^{32}(h_2)&\equiv&\ind\{Y>\zeta\}\left\{
\hat{\Xi}_{\theta}^{(0)}(\tau)R_{\zeta,\psi}^{\tau}(ZZ_2'h_2)
+\hat{\Xi}_{\theta}^{(1)}(\tau)R_{\zeta,\psi}^{\tau}(Z)
R_{\zeta,\psi}^{\tau}(Z_2'h_2)\right\},\\
\hat{\sigma}_{\theta}^{33}(h_3)&\equiv& 
\hat{\Xi}_{\theta}^{(0)}(\tau)R_{\zeta,\psi}^{\tau}(ZZ'h_3)
+\hat{\Xi}_{\theta}^{(1)}(\tau)R_{\zeta,\psi}^{\tau}(Z)
R_{\zeta,\psi}^{\tau}(Z'h_3),\\
\hat{\sigma}_{\theta}^{34}(h_4)&\equiv& 
\hat{\Xi}_{\theta}^{(0)}(\tau)R_{\zeta,\psi}^{\tau}(Zh_4)
+\hat{\Xi}_{\theta}^{(1)}(\tau)R_{\zeta,\psi}^{\tau}(Z)
R_{\zeta,\psi}^{\tau}(h_4),
\end{eqnarray*}
\begin{eqnarray*}
\hat{\sigma}_{\theta}^{41}(h_1)(u)&\equiv&\ind\{Y>\zeta\}
\tilde{Y}(u)e^{r_{\xi}(u;Z,Y)}\left\{
\hat{\Xi}_{\theta}^{(0)}(\tau)h_1+
\hat{\Xi}_{\theta}^{(1)}(\tau)R_{\zeta,\psi}^{\tau}(h_1)\right\},\\
\hat{\sigma}_{\theta}^{42}(h_2)(u)&\equiv&\ind\{Y>\zeta\}
\tilde{Y}(u)e^{r_{\xi}(u;Z,Y)}\left\{
\hat{\Xi}_{\theta}^{(0)}(\tau)Z_2'(u)h_2+
\hat{\Xi}_{\theta}^{(1)}(\tau)R_{\zeta,\psi}^{\tau}(Z_2'h_2)\right\},\\
\hat{\sigma}_{\theta}^{43}(h_3)(u)&\equiv&
\tilde{Y}(u)e^{r_{\xi}(u;Z,Y)}\left\{
\hat{\Xi}_{\theta}^{(0)}(\tau)Z'(u)h_3+
\hat{\Xi}_{\theta}^{(1)}(\tau)R_{\zeta,\psi}^{\tau}(Z'h_3)\right\},\\
\hat{\sigma}_{\theta}^{44}(h_4)(u)&\equiv&
\tilde{Y}(u)e^{r_{\xi}(u;Z,Y)}\left\{
\hat{\Xi}_{\theta}^{(0)}(\tau)h_4(u)+
\hat{\Xi}_{\theta}^{(1)}(\tau)R_{\zeta,\psi}^{\tau}(h_4)\right\},
\end{eqnarray*}
and where
\begin{eqnarray*}
\hat{\Xi}_{\theta}^{(1)}(\tau)&\equiv&\ddot{G}(H^{\theta}(V\wedge\tau))
-\delta\left[\frac{\dddot{G}(H^{\theta}(V\wedge\tau))}
{\dot{G}(H^{\theta}(V\wedge\tau))}-\left\{
\frac{\ddot{G}(H^{\theta}(V\wedge\tau))}
{\dot{G}(H^{\theta}(V\wedge\tau))}\right\}^2\right].
\end{eqnarray*}
Note that all of the above operators are clearly bounded
whenever $\theta$ is bounded. 

The following lemma strengthens
the above G\^{a}teaux derivative to a Fr\'{e}chet derivative.  We will
need this strong differentiability to obtain weak convergence of
our estimators.
\begin{lemma}\label{l3}
Under the regularity conditions of section~2 and for any $\zeta\in[a,b]$
and $\psi_1\in\Psi$, the operator 
$\psi\mapsto P U_{\zeta}^{\tau}(\psi)$
is Fr\'{e}chet differentiable at $\psi_1$,
with derivative $-\psi(\sigma_{\psi_1}(h))$, where
$h$ ranges over ${\cal H}_r$ and is the index
for $P_{\zeta}^{\tau}(\psi)(\cdot)$, $\psi$ ranges over the linear span
$\mbox{lin}\,\Psi$ of $\Psi$, and $0<r<\infty$. 
\end{lemma}

The following lemma gives
us the desired continuous invertibility of both
$\sigma_{\theta_0}$ and the operator 
$\psi\mapsto\psi(\sigma_{\theta_0}(\cdot))$. This last operator will
be needed for weak convergence of regular parameters.
\begin{lemma}\label{l4}
Under the regularity conditions of section~2, 
the linear operator $\sigma_{\theta_0}: \mathcal{H}_{\infty}
\rightarrow \mathcal{H}_{\infty}$ is 
continuously invertible and onto, with inverse $\sigma_{\theta_0}^{-1}$. 
Moreover, the linear operator $\psi\mapsto \psi(\sigma_{\theta_0}(\cdot))$, 
as a map from and to $\mbox{lin}\,\Psi$, 
is also continuously invertible and onto, with inverse
$\psi\mapsto\psi(\sigma_{\theta_0}^{-1}(\cdot))$.
\end{lemma}

\section{The convergence rates of the estimators}

To determine the convergence rates of the estimators, we need to
study closely the log-likelihood process $\tilde{L}_n(\theta)$
near its maximizer. In the parametric setting, this process
can be approximated by its expectation which can be shown to
be locally concave.  For the Cox model, as in \cite{p03}, this
same procedure can be applied to the partial likelihood which
shares the local concavity features of a parametric likelihood.
Unfortunately, in our present set-up, studying the expectation
of $\tilde{L}_n(\theta)$ will lead to problems since $A_0$ has
a density and thus $\Delta A_0(t)=0$ for all $t\in[0,\tau]$.
Hence $\tilde{L}_n(\theta_0)=-\infty$, and
a new approach is needed. The approach we take involves a
careful reparameterization of $\hat{A}_n$. 

From section~4, we know that the maximizer $\hat{A}_n(t)=
\int_0^t\left\{\pp_n W(s;\hat{\theta}_n)\right\}^{-1}$
$\times d\tilde{G}_n(s)$,
where $\tilde{G}_n(t)\equiv\pp_n N(t)$ and $W(\cdot;\cdot)$ is
as defined in~(\ref{c4:e2}). It is easy to see that for all $n$
large enough and all $\theta$ sufficiently close to $\theta_0$,
$t\mapsto\pp_n W(t;\theta)$ is bounded below and above and in
total variation, with large probability.  
Thus, if we use the reparameterization
$\Gamma(\cdot)\mapsto A^{(\Gamma)}_n(\cdot)\equiv\int_0^{(\cdot)}
\exp\{-\Gamma(s)\}d\tilde{G}_n(s)$, and 
maximize $\tilde{L}_n(\xi,A^{(\Gamma)}_n)$ over $\xi$ and $\Gamma$, where 
$\Gamma\in BV$, we will achieve the same NPMLE as before. Note that
the $\Gamma$ component of the maximizer of $\tilde{L}(\xi,A^{(\Gamma)}_n)$
is therefore just $\hat{\Gamma}_n(\cdot)\equiv-\log\pp_n W(\cdot;
\hat{\theta}_n)$.

Define $\Gamma_0(\cdot)\equiv-\log(PW(\cdot;\theta_0))$ and
$\theta_n(\zeta,\gamma,\Gamma)\equiv(\zeta,\gamma,A^{(\Gamma)}_n)$, 
and note that
the reparameterized NPMLE $(\hat{\zeta}_n,\hat{\gamma}_n,\hat{\Gamma}_n)$
is the maximizer of the process
\begin{eqnarray*}
\lefteqn{
(\zeta,\gamma,\Gamma)\mapsto\tilde{X}_n(\zeta,\gamma,\Gamma)
\;\equiv\;\tilde{L}_n(\zeta,\gamma,A^{(\Gamma)}_n)-\tilde{L}_n(\zeta_0,
\gamma_0,A^{(\Gamma_0)}_n)}&&\\
&&\mbox{\hspace{-0.1in}}=\pp_n\left\{\int_0^{\tau}\left[-\Gamma(t)+\Gamma_0(t)
+\log\frac{\dot{G}(H^{\theta_n(\zeta,\gamma,\Gamma)}(t))}
{\dot{G}(H^{\theta_n(\zeta_0,\gamma_0,\Gamma_0)}(t))} +
(r_{\xi}-r_{\xi_0})(t;Z,Y)\right]\right.\\
&&\mbox{\hspace{0.2in}}
\left.\rule[-0.3cm]{0cm}{1.0cm}\times dN(t)
-(G(H^{\theta_n(\zeta,\gamma,\Gamma)}(V))
-G(H^{\theta_n(\zeta_0,\gamma_0,\Gamma_0)}(V)))
\right\}.
\end{eqnarray*}
We will argue shortly that $\tilde{X}_n$ is uniformly consistent for
the function
\begin{eqnarray*}
\lefteqn{(\zeta,\gamma,\Gamma)\mapsto\tilde{X}(\zeta,\gamma,\Gamma)}&&\\
&\equiv&P\left\{\int_0^{\tau}\left[-\Gamma(t)+\Gamma_0(t)
+\log\frac{\dot{G}(H^{\theta_0(\zeta,\gamma,\Gamma)}(t))}
{\dot{G}(H^{\theta_0}(t))} +
(r_{\xi}-r_{\xi_0})(t;Z,Y)\right]dN(t)\right.\\
&&\left.\rule[-0.3cm]{0cm}{1.0cm}-(G(H^{\theta_0(\zeta,\gamma,\Gamma)}(V))
-G(H^{\theta_0}(V)))
\right\},
\end{eqnarray*}
where $\theta_0(\zeta,\gamma,\Gamma)\equiv(\zeta,\gamma,
A^{(\Gamma)}_0)$, $A^{(\Gamma)}_0(\cdot)\equiv\int_0^{(\cdot)}
\exp\{-\Gamma(s)\}d\tilde{G}_0(s)$, and $\tilde{G}_0(t)\equiv P N(t)$.
It will occasionally be useful to use the shorthand
$\lambda\equiv(\gamma,\Gamma)$, 
$\hat{\lambda}_n\equiv(\hat{\gamma}_n,\hat{\Gamma}_n)$ and 
$\lambda_0\equiv(\gamma_0,\Gamma_0)$.

Define the modified parameter space 
$\Theta^{\ast}\equiv (a,b)\times\Upsilon\times B_2\times B_1
\times BV$; and, for each
$h=(h_1,h_2,h_3,h_4,h_5)\in\re\times{\cal H}_{\infty}$, define
the metric $\rho_2(h)\equiv(|h_1|+|h_2|^2+\|h_3\|^2+\|h_4\|^2
+\|h_5\|_{\infty}^2)^{1/2}$, where $\|\cdot\|_{\infty}$ is the uniform
norm. Note that $|h_1|$ is deliberately not squared.
For each $\epsilon>0$ and $k<\infty$, define $B_{\epsilon}^{\ast k}
\equiv\{(\zeta,\lambda)\in\Theta^{\ast}:
\rho_2((\zeta,\lambda)-(\zeta_0,\lambda_0))<\epsilon,\|\Gamma\|_v\leq k\}$.
Note that for some $k_0<\infty$ and
any $\epsilon>0$, $(\hat{\zeta}_n,\hat{\lambda}_n)$ 
is eventually in $B_{\epsilon}^{\ast k_0}$ for all 
$n$ large enough by theorem~\ref{t1} above combined with 
lemma~\ref{l5} below:
\begin{lemma}\label{l5}
There exists a $k_0<\infty$ such that 
$\limsup_{n\rightarrow\infty}\|\hat{\Gamma}_n\|_v\leq k_0$ and
$\lim_{n\rightarrow\infty}\|\hat{\Gamma}_n-\Gamma_0\|_{\infty}=0$ 
outer almost surely.
\end{lemma}

Now we study the local behavior of $\tilde{X}$. First
fix $\zeta\in(a,b)$. Since, for any $g\in BV$, 
\[\left.\frac{\partial A^{(\Gamma+t g)}_0(\cdot)}{\partial t}\right|_{t=0}
=-\int_0^{(\cdot)}g(s)dA^{(\Gamma)}_0(s),\]
we obtain that the first derivative of $(\gamma,\Gamma)\mapsto
\tilde{X}(\zeta,\gamma,\Gamma)$ in the direction $h\in{\cal H}_{\infty}$,
is precisely $-PU_{\zeta}^{\tau}(\gamma,A^{(\Gamma)}_0)(h)$. Moreover, 
by definition of the score and information operators,
the second derivative in the same direction is 
$-\psi^h_{\Gamma}\left(\sigma_{\left(\zeta,\gamma,
A^{(\Gamma)}_0\right)}(h)\right)$, 
where $\psi^h_{\Gamma}\equiv\left(h_1,h_2,h_3,\int_0^{(\cdot)}
h_4(s)dA^{(\Gamma)}_0(s)\right)$. At the point $(\zeta,\gamma,\Gamma)=
(\zeta_0,\gamma_0,\Gamma_0)$, the first derivative is $0$, 
while the second derivative is $<0$, by lemma~\ref{l4}.
By the smoothness of the score
and information operators ensured by condition~D1 and~D2, 
and by the arbitrariness
of $h$, we now have that the function 
$(\gamma,\Gamma)\mapsto\tilde{X}(\zeta,\gamma,\Gamma)$
is concave for every $(\zeta,\gamma,\Gamma)\in B_{\epsilon}^{\ast k_0}$, for
sufficiently small $\epsilon$.

Now note that $\tilde{X}(\zeta,\gamma,\Gamma)=P l^{\ast}(\zeta,\gamma,
\Gamma)-P l^{\ast}(\zeta_0,\gamma_0,\Gamma_0)$, where 
$l^{\ast}(\zeta,\gamma,\Gamma)\equiv$
\begin{eqnarray}
&&\;\;-\int_0^{\tau}\Gamma(t)dN(t)
+l_1^{\psi(\gamma,\Gamma)}(V,\delta,Z)\ind\{Y\leq\zeta\}
+l_2^{\psi(\gamma,\Gamma)}(V,\delta,Z)
\ind\{Y>\zeta\},\label{new.j14.e1}
\end{eqnarray}
and where $l_j^{\psi}$, $j=1,2$, are as defined in section~3, and
$\psi(\gamma,\Gamma)\equiv(\gamma,A^{(\Gamma)}_0)$.
By condition~B2, we now have that for small enough $\epsilon>0$,
$\zeta\mapsto\tilde{X}(\zeta,\gamma,\Gamma)$ 
is right and left continuously differentiable for all 
$(\zeta,\gamma,\Gamma)\in B_{\epsilon}^{\ast k_0}$,
with left partial derivative 
\[\dot{X}_{\zeta}^{-}(\gamma,\Gamma)
\equiv P\left\{\left.l_1^{\psi(\gamma,\Gamma)}(V,\delta,Z)
-l_2^{\psi(\gamma,\Gamma)}(V,\delta,Z)\right|Y=\zeta\right\}\]
and right partial derivative 
\[\dot{X}_{\zeta}^{+}(\gamma,\Gamma)
\equiv P\left\{\left.l_1^{\psi(\gamma,\Gamma)}(V,\delta,Z)
-l_2^{\psi(\gamma,\Gamma)}(V,\delta,Z)\right|Y=\zeta+\right\}.\]

We now have the following lemmas on the local behavior of $\tilde{X}$
with respect to $\zeta$:
\begin{lemma}\label{l6}
Under the conditions of section~2, 
$\dot{X}_{\zeta_0}^{-}(\gamma_0,\Gamma_0)>0$ and
$\dot{X}_{\zeta_0}^{+}(\gamma_0,\Gamma_0)<0$.
\end{lemma}
\begin{lemma}\label{l7}
There exists $\epsilon_1,k_1>0$ such that
$\tilde{X}(\zeta,\gamma,\Gamma)\leq -k_1|\zeta-\zeta_0|$ 
for all $(\zeta,\gamma,\Gamma)\in B_{\epsilon_1}^{\ast k_0}$.
\end{lemma}

The two previous lemmas can be combined with the next lemma, 
lemma~\ref{l8}, to yield $\sqrt{n}$ rates for all of the parameters 
(theorem~\ref{t.l9}):
\begin{lemma}\label{l8}
There exists an $\epsilon_2>0$ such that $D_n\equiv\sqrt{n}(\tilde{X}_n
-\tilde{X})$ converges weakly to a tight mean zero Gaussian process
$D_0$, in $\ell^{\infty}(B_{\epsilon_2}^{\ast k_0})$, for which
$D_0(\zeta,\gamma,\Gamma)\rightarrow 0$ in probability, as
$\rho_2((\zeta,\gamma,\Gamma)-(\zeta_0,\gamma_0,\Gamma_0))
\rightarrow 0$.$\Box$
\end{lemma}

\begin{theorem}\label{t.l9}
Under the conditions of section~2,
$\sqrt{n}|\hat{\zeta}_n-\zeta_0|=O_P(1)$,
$\sqrt{n}\|\hat{\psi}_n-\psi_0\|_{\infty}=O_P(1)$, and
$\sqrt{n}\|\hat{\Gamma}_n-\Gamma_0\|_{\infty}=O_P(1)$.
\end{theorem}

To refine the rate for $\hat{\zeta}_n$, we need two more lemmas,
lemmas~\ref{l10} and~\ref{l11} below. We will
also need to define the process $\zeta\mapsto\tilde{X}_n^{\ast}(\zeta)
\equiv$
\begin{eqnarray*}
&&\pp_n\left\{\int_0^{\tau}\left[
\log\frac{\dot{G}(H^{\theta_0(\zeta,\gamma_0,\Gamma_0)}(t))}
{\dot{G}(H^{\theta_0}(t))} +
(r_{(\zeta,\gamma_0)}-r_{\xi_0})(t;Z,Y)\right]dN(t)\right.\\
&&\mbox{\hspace{0.4in}}
\left.\rule[-0.3cm]{0cm}{1.0cm}-(G(H^{\theta_0(\zeta,\gamma_0,\Gamma_0)}(V))
-G(H^{\theta_0}(V)))\right\}.
\end{eqnarray*}
\begin{lemma}\label{l10}
$0\leq \tilde{X}_n(\hat{\zeta}_n,\hat{\lambda}_n)-
\tilde{X}_n^{\ast}(\hat{\zeta}_n)\leq O_P(n^{-1})$.
\end{lemma}
\begin{lemma}\label{l11}
There exists an $\epsilon_3>0$ and $k_2<\infty$ such that,
for all $0\leq\epsilon\leq\epsilon_3$ and $n\geq 1$,
$\Exp{\sup_{|\zeta-\zeta_0|\leq\epsilon}|\tilde{D}_n(\zeta)|}
\leq k_2\sqrt{\epsilon}$,
where $\tilde{D}_n(\zeta)\equiv\sqrt{n}(\tilde{X}_n^{\ast}(\zeta)
-\tilde{X}(\zeta,\lambda_0))$.
\end{lemma}

We now have the following theorem about the convergence rate for
$\hat{\zeta}_n$:
\begin{theorem}\label{t2}
Under the conditions of section~2, $n|\hat{\zeta}_n-\zeta_0|=O_P(1)$.
\end{theorem}

{\it Proof.} The method of proof involves a ``peeling device'' (see, 
for example, the proof of theorem~5.1 of \cite{ih81},
or the proof of theorem~2 of \cite{p03}).
Fix $\epsilon>0$. By consistency and lemma~\ref{l5}, 
$P((\hat{\zeta}_n,\hat{\lambda}_n)\in B_{\epsilon_4}^{\ast k_0})
\geq 1-\epsilon$ for
all $n$ large enough, where $\epsilon_4=
\epsilon_1\wedge \epsilon_2\wedge\epsilon_3$.
By lemma~\ref{l10}, there exists an $M_1^{\ast}<\infty$ such that
$P(\tilde{X}_n(\hat{\zeta}_n,\hat{\lambda}_n)-
\tilde{X}_n^{\ast}(\hat{\zeta}_n)>M_1^{\ast}/n)\leq\epsilon$. 
For integers $k\geq 1$, let $m_k\equiv k^4$. We now have, for any
integer $k\geq 1$, that 
$\limsup_{n\rightarrow\infty}P\left(n|\hat{\zeta}_n-\zeta_0|>m_k\right)$
\begin{eqnarray}
&\leq&\limsup_{n\rightarrow\infty}P\left(
n|\hat{\zeta}_n-\zeta_0|>m_k,\;(\hat{\zeta}_n,\hat{\lambda}_n)
\in B_{\epsilon_4}^{\ast k_0},\right.\nonumber\\
&&\left.\tilde{X}_n(\hat{\zeta}_n,\hat{\lambda}_n)
-\tilde{X}_n^{\ast}(\hat{\zeta}_n)
\leq \frac{M_1^{\ast}}{n}\right)+2\epsilon\nonumber\\
&\leq&\limsup_{n\rightarrow\infty}P\left(\sup_{\zeta:\,m_k/n<|\zeta-\zeta_0|
\leq\epsilon_4}\tilde{X}_n^{\ast}(\zeta)\geq -\frac{M_1^{\ast}}{n}\right)
+2\epsilon\nonumber\\
&\leq&\limsup_{n\rightarrow\infty}\sum_{j=k}^{k_{\epsilon_4}}
P\left(\sup_{\zeta:\,m_j/n<|\zeta-\zeta_0|\leq (m_{j+1}/n)
\wedge\epsilon_4}\tilde{D}_n(\zeta)\right.\label{t2.e1}\\
&&\left.\mbox{\hspace{0.4in}}\geq\sqrt{n}\left(\frac{k_1m_j}{n}
-\frac{M_1^{\ast}}{n}\right)\right)+2\epsilon,\nonumber
\end{eqnarray}
by lemma~\ref{l7}, where $k_{\epsilon_4}=
\min\{k:\,m_{k+1}\geq n\epsilon_4\}$. But, by lemma~\ref{l11}, 
\[\mbox{(\ref{t2.e1})}\leq\limsup_{n\rightarrow\infty}
\sum_{j=k}^{k_{\epsilon_4}}\frac{k_2\sqrt{m_{j+1}}}{k_1m_j-M_1^{\ast}}
+2\epsilon\leq\sum_{j=k}^{\infty}\frac{k_2(j+1)^2}
{k_1j^4-M_1^{\ast}}+2\epsilon.\]
We can now choose $k<\infty$ large enough so that this last term
$\leq 3\epsilon$. Since $\epsilon>0$ was arbitrary, we now have that
$\lim_{m\rightarrow\infty}\limsup_{n\rightarrow\infty}
P(n|\hat{\zeta}_n-\zeta_0|>m)=0$, and the desired conclusion follows.$\Box$

\section{Weak convergence of the estimators} 

\subsection{The asymptotic distribution of the change-point
estimator} 

Denote $\mathbb{U}_{n,M}\equiv\{u=n(\zeta-\zeta_0):\zeta\in[a,b],|u|\leq M\}$
and $\zeta_{n,u}\equiv\zeta_0+u/n$.
The limiting distribution of $n(\hat{\zeta}_n-\zeta_0)$
will be deduced from the behavior of the restriction of
the process $u \rightarrow n[\tilde{L}_n(\hat{\psi}_n,\zeta_{n,u})
-\tilde{L}_n(\hat{\psi}_n,\zeta_0)]$ to
the compact set $\mathbb{U}_{n,M}$, for $M$ sufficiently large.
\begin{theorem}\label{t3}
The following approximation holds for all $M > 0$, as $n \rightarrow \infty$:
\[u\mapsto n[\tilde{L}_n(\hat{\psi}_n,\zeta_{n,u}) 
-\tilde{L}_n (\hat{\psi}_n,\zeta_0)]= 
Q_n(u)+o_P^{\mathbb{U}_{n,M}}(1),\]
where $o_P^B(1)$ denotes a term going to zero in probability uniformly
over the set $B$ and $\mbox{Q}_n(u)=$ 
\begin{eqnarray*}
n\mathbb{P}_n
\left\{\left(\ind\{\zeta_{n,u}<Y \le \zeta_0\} 
- \ind\{\zeta_0<Y \le \zeta_{n,u}\}\right)
\left[l_2^{\psi_0}(V,\delta,Z)-l_1^{\psi_0}(V,\delta,Z)\right]\right\}.
\end{eqnarray*}
\end{theorem}

Let $Q_n(u)=Q_n^+(u)\ind\{u>0\}-Q_n^-(u)\ind\{u<0\}$. 
We now study the weak convergence of $Q_n$ as a random variable on the
space of cadlag functions $D$ with the Skorohod topology, and on
its restriction to the space $D_M$ of cadlag functions on $[-M,
M]$, for any $M > 0$, similar to the approach taken in \cite{p03}. 
In order to describe the asymptotic distribution of $Q_n$,
let $\nu^+$ and $\nu^-$ be two independent jump processes on~$\mathbb{R}$
such that $\nu^+ (s)$ is a Poisson variable with parameter 
$s^+\tilde{h}(\zeta_0)$ and $\nu^- (s)$ is a Poisson variable with
parameter $(-s)^+\tilde{h}(\zeta_0)$. Here, 
$u^+$ denotes $u\vee 0$. Let $(\check{V}_k^+)_{k\ge 1}$ 
and $(\check{V}_k^-)_{k \ge 1}$ be independent sequences of i.i.d.
random variables with characteristic functions
\[\phi^+(t)=P\left[e^{it\check{V}_k^+}\right]
=P\left[\left.e^{it\left\{l_1^{\psi_0}(V,\delta,Z)-l_2^{\psi_0}(V,\delta,Z)
\right\}}\right|Y=\zeta_0^+\right],\]
and
\[\phi^-(t)=P\left[e^{it\check{V}_k^-}\right]
=P\left[\left.e^{it\left\{l_1^{\psi_0}(V,\delta,Z)-l_2^{\psi_0}(V,\delta,Z)
\right\}}\right|Y=\zeta_0\right],\]
respectively, where $(\check{V}_k^+)_{k \ge 1}$ and 
$(\check{V}_k^-)_{k \ge 1}$ are independent of $\nu^+$ and $\nu^-$.

Let $Q(s) = Q^+(s)\ind\{s > 0\} - Q^-(s)\ind\{s < 0\}$ be the
right-continuous jump process defined by
\[Q^+(s)=\sum_{0 \le k \le {\nu}^+(s)}\check{V}_k^+, ~ ~ ~ 
Q^-(s)=\sum_{0 \le k \le {\nu}^-(s+)}\check{V}_k^-, \]
where $\check{V}_0^+=\check{V}_0^-=0$. 
Using a modification of the arguments in \cite{p03}, we obtain:
\begin{theorem}\label{t4}
Under the regularity conditions of section~2,
the process $Q_n$ converges weakly to $Q$ in $D_M$, for every $M > 0$; 
$n(\hat{\zeta}_n - \zeta_0) = \argmax_{u}Q_n(u) + o_p(1)$ which
converges weakly to $\hat{v}_Q \equiv \argmin\{|v| : Q(v) =
\argmax\,Q\}$; and $n(\hat{\zeta}_n-\zeta_0)$ and
$\sqrt{n}\pp_nU_{\zeta_0}^{\tau}(\psi_0)(h)$ are asymptotically
independent for all $h\in{\cal H}_{\infty}$.
\end{theorem}

\subsection{Asymptotic normality of the regular parameters} We use
Hoffmann-J{\o}rgensen weak convergence as described in \cite{vw96}. 
We have the
following result:
\begin{theorem}\label{t5}
Under the conditions of theorem 1,  $\sqrt{n}(\hat \psi_n -
\psi_0)$ is asymptotically linear, with influence function $\tilde
l(h) = U_{\zeta_0}^{\tau}(\psi_0)(\sigma_{\theta_0}^{-1}(h))$, $h
\in{\cal H}_1$, converging weakly in the uniform norm to a tight, mean
zero Gaussian process $\mathbb{Z}$ with covariance
$E[\tilde l(g) \tilde l(h)]$, for all $g, h \in H_1$. Thus
$n(\hat{\zeta}_n-\zeta_0)$ and $\sqrt{n}(\hat{\psi}_n-\psi_0)$
are asymptotically independent.
\end{theorem}

\begin{remark}\label{r1}
Since $\sqrt{n}(\hat \psi_n - \psi_0)$ is asymptotically linear,
with influence function contained in
the closed linear span of the tangent space (since
$\sigma_{\theta_0}$ is continuously invertible), $\hat\psi_n$ is
regular and hence as efficient as if $\zeta_0$ were known, by
Theorem 5.2.3 and Theorem 5.2.1 of \cite{bkrw98}.
\end{remark}

\section{Inference when $\alpha_0\neq 0$ or $\eta_0\neq 0$}
In this section we develop Monte Carlo methods for inference for the parameter estimators when
it is known that either $\alpha_0\neq 0$ or $\eta_0\neq 0$, i.e., it is known that condition~C2
is satisfied. In section~9, 
we develop a hypothesis testing procedure to assess whether
$H_0:\alpha_0=0=\eta_0$ holds (i.e., that~C2 does not hold). When it is known that $H_0$ holds, 
the model reduces to the usual transformation model 
(see \cite{sv04}), 
and thus validity of the bootstrap will follow from arguments
similar to those used in the proof of 
corollary~1 of \cite{klf04}. 

\subsection{Inference for the change-point} One possibility for
inference for $\zeta$ is to use the subsampling bootstrap \cite{pr94}
which is guaranteed to work, provided the subsample
sizes $\ell_n$ satisfy $\ell_n\rightarrow\infty$ and $\ell_n/n\rightarrow 0$.
However, this approach is very computationally intense since,
for each subsample, the likelihood must be maximized over the entire
parameter space. To ameliorate the computational strain, we propose as
an alternative the following specialized parametric bootstrap. 
Let $\tilde{F}_+$ and $\tilde{F}_-$ be the 
distribution functions corresponding to the moment generating functions
$\phi^+$ and $\phi^-$, respectively. We need to make the following
additional assumption:
\begin{enumerate}
\item[B5:] Both $\tilde{F}_+$ and $\tilde{F}_-$ are continuous.
\end{enumerate}
Now let $\tilde{m}_n$ be the minimum of the number of $Y$ observations
in the sample $>\hat{\zeta}_n$ and the number of $Y$ observations $<\hat{\zeta}_n$. Now
choose sequences of possibly data dependent integers $1\leq C_{1,n}<C_{2,n}\leq \tilde{m}_n$ 
such that $C_{1,n}\rightarrow\infty$, 
$C_{2,n}-C_{1,n}\rightarrow\infty$, and 
$C_{2,n}/n\rightarrow 0$, in probability, as $n\rightarrow\infty$. 
Note that if one
chooses $C_{1,n}$ to be the closest integer to $\tilde{m}_n^{1/4}$ and $C_{2,n}$ to be
the closest integer to $\tilde{m}_n^{3/4}$, 
the given requirements will be satisfied since 
$\tilde{m}_n\rightarrow\infty$, in probability, by assumption~B1. Let
$X_{(1)},\ldots,X_{(n)}$ be the complete data observations corresponding
to the order statistics $Y_{(1)},\ldots,Y_{(n)}$ of the $Y$ observations. 
Also let $\tilde{k}_n\equiv C_{2,n}-C_{1,n}+1$, and define $\tilde{l}_n$ 
to be the integer satisfying $\hat{\zeta}_n=Y_{(\tilde{l}_n)}$. 
The existence of this integer follows from the form of the MLE. 

Now, for $j=1,\ldots,\tilde{k}_n$, and any $\psi\in\Psi$, define 
\begin{eqnarray*}
\check{V}_{j,\psi}^+&\equiv& l_1^{\psi}(
V_{(\tilde{l}_n+C_{1,n}+j-1)},\delta_{(\tilde{l}_n+C_{1,n}+j-1)},Z_{(\tilde{l}_n+C_{1,n}+j-1)})\\
&&-l_2^{\psi}(V_{(\tilde{l}_n+C_{1,n}+j-1)},\delta_{(\tilde{l}_n+C_{1,n}+j-1)},
Z_{(\tilde{l}_n+C_{1,n}+j-1)}),\\
\check{V}_{j,\psi}^-&\equiv&l_1^{\psi}(
V_{(\tilde{l}_n-C_{1,n}-j)},\delta_{(\tilde{l}_n-C_{1,n}-j)},Z_{(\tilde{l}_n-C_{1,n}-j)})\\
&&-l_2^{\psi}(V_{(\tilde{l}_n-C_{1,n}-j)},\delta_{(\tilde{l}_n-C_{1,n}-j)},
Z_{(\tilde{l}_n-C_{1,n}-j)}),
\end{eqnarray*}
$Y^+_j\equiv Y_{(\tilde{l}_n+C_{1,n}+j-1)}$, and $Y^-_j\equiv Y_{(\tilde{l}_n-C_{1,n}-j)}$. Also let $\hat{F}_+^n$ be 
the data-dependent distribution function
for a random variable drawn with replacement from 
$\{\check{V}_{1,\hat{\psi}_n}^+,\ldots,\check{V}_{\tilde{k}_n,
\hat{\psi}_n}^+\}$, and let $\hat{F}_-^n$ be 
the data-dependent distribution function
for a random variable drawn with replacement from 
$\{\check{V}_{1,\hat{\psi}_n}^-,\ldots,$ $\check{V}_{\tilde{k}_n,
\hat{\psi}_n}^-\}$. By the smoothness 
of the terms involved, it is easy to verify
that both $\sup_{1\leq j\leq\tilde{k}_n}$ $\left|\check{V}_{j,\hat{\psi}_n}^+
-\check{V}_{j,\psi_0}^+\right|=o_P(1)$ and 
$\sup_{1\leq j\leq\tilde{k}_n}\left|\check{V}_{j,\hat{\psi}_n}^-
-\check{V}_{j,\psi_0}^-\right|=o_P(1)$. Moreover, by assumption~B2(i), the
fact that $n(\hat{\zeta}_n-\zeta_0)=O_P(1)$, and the conditions on $C_{1,n}$
and $C_{2,n}$, we have that both $P(Y^-_{1}<\zeta_0<Y^+_{1})\rightarrow 1$ and
$Y^+_{\tilde{k}_n}-Y^-_{\tilde{k}_n}=o_P(1)$. Thus, by assumption~B2(ii),
the collection $\{\check{V}^+_{1,\psi_0},\ldots,$
$\check{V}^+_{\tilde{k}_n,\psi_0}\}$ converges
in distribution to an i.i.d. sample of random 
variables with characteristic function
$\phi^+$, while the collection $\{\check{V}^-_{1,\psi_0},\ldots,
\check{V}^-_{\tilde{k}_n,\psi_0}\}$ is 
independent of the first collection and converges
in distribution to an i.i.d. sample of random 
variables with characteristic function $\phi^-$.  By assumption~B5
and the fact that $\tilde{k}_n\rightarrow\infty$, in probability, 
we now have that both
$\sup_{v\in\re}|\hat{F}_+^n(v)-\tilde{F}_+(v)|=o_P(1)$ and
$\sup_{v\in\re}|\hat{F}_-^n(v)-\tilde{F}_-(v)|=o_P(1)$.

Now let $\hat{h}_n$ be a consistent estimator of $\tilde{h}(\zeta_0)$.
Such an estimator can be obtained from a kernel density estimator of
$\tilde{h}$ based on the $Y$ observations and evaluated at $\hat{\zeta}_n$.
The basic idea of our parametric bootstrap is to create a stochastic
process $\hat{Q}_n$ defined similarly to the process $Q$ described
in section~7.1. To this end,
let $\hat{\nu}^+$ and $\hat{\nu}^-$ be two independent jump processes
defined on the interval $\tilde{B}_n\equiv
[-n(\hat{\zeta}_n-a),n(b-\hat{\zeta}_n)]$ 
such that $\hat{\nu}^+(s)$ is Poisson with parameter $s^+\hat{h}_n$
and $\hat{\nu}^-(s)$ is Poisson with parameter $(-s)^+\hat{h}_n$.
Also let $(\check{V}_{\ast,k}^+)_{k\geq 1}$ and 
$(\check{V}_{\ast,k}^-)_{k\geq 1}$ be two independent sequences of
i.i.d. random variables drawn from $\hat{F}_+^n$ and $\hat{F}_-^n$
and independent of the Poisson processes. Now construct 
$u\mapsto\hat{Q}_n(u)\equiv\hat{Q}_n^+(u)\ind\{u>0\}
-\hat{Q}_n^-(u)\ind\{u<0\}$ on the interval $\tilde{B}_n$, 
where $\hat{Q}_n^+(u)\equiv
\sum_{0\leq k\leq\hat{\nu}^+(u)}\check{V}_{\ast,k}^+$ and
$\hat{Q}_n^-(u)\equiv\sum_{0\leq k\leq\hat{\nu}^-(u+)}\check{V}_{\ast,k}^-$.
Finally, we compute $\hat{v}_{\ast}\equiv\argmin_{\tilde{B}_n}\left\{|v|:
\hat{Q}_n(v)=\argmax_{\tilde{B}_n}\hat{Q}_n\right\}$.
The following proposition now follows from 
the fact that $P(K\in\tilde{B}_n)\rightarrow 1$ for all
compact $K\subset\re$:
\begin{proposition}\label{p1}
The conditional distribution of $\hat{v}_{\ast}$ given the data is
asymptotically equal to the distribution of $\hat{v}_Q$ defined
in theorem~\ref{t4}. 
\end{proposition}

Hence for any $\pi>0$, we can consistently estimate the 
$\pi/2$ and $1-\pi/2$ quantiles of $\hat{v}_Q$ based
on a large number of independent draws from $\hat{v}_{\ast}$,
which estimates we will denote by $\hat{q}_{\pi/2}$ and
$\hat{q}_{1-\pi/2}$, respectively. Thus an asymptotically
valid $1-\pi$ confidence interval for $\zeta_0$ is
$[\hat{\zeta}_n-\hat{q}_{1-\pi/2},\hat{\zeta}_n-\hat{q}_{\pi/2}]$.

\subsection{Inference for regular parameters} Because $\hat{\zeta}_n$
is $n$-consistent for $\zeta_0$, $\zeta_0$ can be treated as known
in constructing inference for the regular parameters. 
Accordingly, we propose bootstrapping the likelihood and maximizing
over $\psi$ while holding $\zeta$ fixed at $\hat{\zeta}_n$. This will
significantly reduce the computational demands of the bootstrap. 
Also, to avoid the occurrence of ties during resampling,
we suggest the following weighted bootstrap alternative
to the usual nonparametric bootstrap. First generate
$n$ i.i.d. positive random variables $\kappa_1,\ldots,\kappa_n$,
with mean $0<\mu_{\kappa}<\infty$, variance
$0<\sigma_{\kappa}^2<\infty$, and with
$\int_0^{\infty}\sqrt{P(\kappa_1>u)}du<\infty$. Divide each weight
by the sample average of the weights $\bar{\kappa}$, to obtain
``standardized weights'' $\kappa_1^{\circ},\ldots,\kappa_n^{\circ}$
which sum to~$n$. For a real, measurable function $f$, define the 
weighted empirical measure $\pp_n^{\circ}f\equiv n^{-1}
\sum_{i=1}^n\kappa_i^{\circ}f(X_i)$. Recall that the nonparametric bootstrap
empirical measure $\pp_n^{\bullet}f\equiv n^{-1}\sum_{i=1}^n
\kappa_i^{\bullet}f(X_i)$ uses multinomial weights
$\kappa_1^{\bullet},\ldots,\kappa_n^{\bullet}$, 
where $\Exp{\kappa_i^{\bullet}}=1$, $i=1,\ldots,n$, and
$\sum_{i=1}^n\kappa_i^{\bullet}=n$ almost surely.

The proposed weighted bootstrap estimate $\hat{\psi}_n^{\circ}$
is obtained by maximizing $\tilde{L}_n^{\circ}(\psi,\hat{\zeta}_n)$ over
$\psi\in\Psi$, where $\tilde{L}_n^{\circ}$ is obtained by replacing
$\pp_n$ with $\pp_n^{\circ}$ in the definition of $\tilde{L}_n$ 
from section~3. We can similarly defined a modified nonparametric
bootstrap $\hat{\psi}_n^{\bullet}$ as the $\argmax$ of
$\psi\mapsto\tilde{L}_n^{\bullet}(\psi,\hat{\zeta}_n)$, where
$\tilde{L}_n^{\bullet}$ is obtained by replacing $\pp_n$ with
$\pp_n^{\bullet}$ in the definition of $\tilde{L}_n$. The following
corollary establishes the validity of both kinds 
of bootstraps:
\begin{corollary}\label{c1}
Under the conditions of theorem~\ref{t5}, the conditional
bootstrap of $\hat{\psi}_n$, based on either 
$\hat{\psi}_n^{\bullet}$ or $\hat{\psi}_n^{\circ}$, 
is asymptotically consistent for
the limiting distribution $\mathbb{Z}$ in the following sense:
Both $\sqrt{n}(\hat{\psi}_n^{\bullet}-\hat{\psi}_n)$ and
$\sqrt{n}(\mu_{\kappa}/\sigma_{\kappa})(\hat{\psi}_n^{\circ}
-\hat{\psi}_n)$ are asymptotically measurable, and both
\begin{enumerate}
\item[(i)] $\sup_{g\in BL_1}\left|E_{\bullet}g\left(
\sqrt{n}(\hat{\psi}_n^{\bullet}-\hat{\psi}_n)\right)
-Eg(\mathbb{Z})\right|\rightarrow 0$ in outer probability and
\item[(ii)] $\sup_{g\in BL_1}\left|E_{\circ}g\left(
\sqrt{n}(\mu_{\kappa}/\sigma_{\kappa})
(\hat{\psi}_n^{\circ}-\hat{\psi}_n)\right)
-Eg(\mathbb{Z})\right|\rightarrow 0$ in outer probability, 
\end{enumerate}
where $BL_1$ is the space of functions mapping
$\re^{d+q+1}\times\ell^{\infty}[0,\tau]\mapsto\re$ which are 
bounded in absolute value by~1 and have Lipschitz norm $\leq 1$.
Here, $E_{\bullet}$ and $E_{\circ}$ are expectations that are
taken over the multinomial and standardized weights, respectively,
conditional on the data.
\end{corollary}

\begin{remark}\label{r2}
As discussed in remark~15 of \cite{klf04}, the
choice of weights $\kappa_1,\ldots,\kappa_n$ in this kind of
setting does not effect the first order asymptotics. However,
it may have an effect on finite samples. In our experience, we
have found that both exponential and truncated exponential weights
perform quite well.
\end{remark}

\section{Test for the presence of a change-point}
Constructing a valid
test of the null hypothesis that there is no change-point,
$H_0:\alpha_0=0=\eta_0$, poses an interesting challenge.
Since the location of the change-point is no longer identifiable
under $H_0$, this is an example of the issue studied in 
\cite{a01}. The test statistic we propose is
a functional of the $\alpha$ and $\eta$ components of the score process, 
$\zeta\mapsto 
\hat{S}_{1}(\zeta)\equiv\sqrt{n}\pp_n
(U^{\tau}_{\zeta, 1}(\hat \psi_0), 
U^{\tau}_{\zeta,2}(\hat \psi_0)')'$, where $\zeta\in[a,b]$,
$\hat\psi_0 \equiv (0,0, \hat\beta_0,\hat A_0)$, 
and where $(\hat{\beta}_0,\hat{A}_0)$ 
is the restricted MLE of $(\beta_0, A_0)$ under the 
assumption that $\alpha=0$ and $\eta=0$. This MLE is relatively easy to
compute since estimation of $\zeta$ is not needed. Specifically,
we have from section~3, that $\hat{\psi}_0$ is the maximizer of 
\begin{eqnarray}
\psi&\mapsto&\pp_n\left\{\delta\log(n\Delta A(V))+l_1^{\psi}(V,\delta,Z)
\right\}.\label{s9.e1}
\end{eqnarray}
We also define for future use  
$h\mapsto\hat{S}_{2}(h)\equiv\sqrt{n}\pp_n
(U^{\tau}_{\zeta,3}(\hat\psi_0)(h_3),U^{\tau}_{\zeta,4}(\hat\psi_0)(h_4))'$, 
where $h\in{\cal H}_1$. The statistic we propose using is $\hat T_n\equiv\sup
_{\zeta \in [a, b]}\left\{\hat{S}_{1}'(\zeta)\hat V_n^{-1}(\zeta)\right.$
$\left.\times\hat{S}_{1} (\zeta)\right\}$, where
$\hat V_n (\zeta)$ is a consistent
estimator of the covariance of $\hat{S}_{1}(\zeta)$.

There are several reasons for us to consider the sup functional of score 
statistics instead of wald or likelihood ratio statistics. Firstly, the score 
statistic is much less computational intense which makes the bootstrap 
implementation feasible. Secondly, we choose the sup functional because of
its guarantee to have some power under local alternatives, as argued in 
\cite{d87} and which we prove below. We note, however, that \cite{ap94}
argue that certain weighted averages of score statistics are optimal 
tests in some settings. A careful analysis of the relative merits of
the two approaches in our setting is beyond the scope of the current paper 
but is an interesting topic for future research. However, as a step in
this direction, we will compare $\hat{T}_n$ with the integrated statistic
$\tilde{T}_n\equiv\int_{[a,b]}\left\{\hat{S}_{1} ' (\zeta) 
\hat V_n^{-1}(\zeta) \hat{S}_{1} (\zeta)\right\}d\zeta$.

In this section, we first discuss a Monte Carlo technique which 
enables computation of $\hat{V}_n(\zeta)$, so that 
$\hat{T}_n$ and $\tilde{T}_n$ can be calculated in the first place, 
as well as computation of critical values for hypothesis testing.
We then discuss the asymptotic properties of the statistics
under a sequence of contiguous alternatives so that power can be
verified. Specifically, we assume that all the conditions of section~2
hold except for C2 which we replace with
\begin{enumerate}
\item[C2':] For each $n\geq 1$, $\alpha_0=\alpha_{\ast}/\sqrt{n}$ and
$\eta_0=\eta_{\ast}/\sqrt{n}$, for some fixed $\alpha_{\ast}\in\re$
and $\eta_{\ast}\in\re^q$. The joint distribution of $(C,Z,Y)$ does
not change with $n$.
\end{enumerate}
Note that when $\alpha_{\ast}\neq 0$ or $\eta_{\ast}\neq 0$,
condition~C2' will cause the distribution of the failure time $T$,
given the covariates $(Z,Y)$, to change with $n$, and the
value of $\zeta_0$ will affect this distribution. 

\subsection{Monte Carlo computation and inference} 
While the nonparametric bootstrap may be a reasonable approach,
it is unclear how to verify its theoretical properties in this context.
We will use instead the weighted bootstrap, based on the multipliers
$\kappa_1^{\circ},\ldots,\kappa_n^{\circ}$ defined in section~8.2.
Let $\pp_n^{\circ}$ be the corresponding weighted empirical measure,
and define $\hat{\psi}_0^{\circ}$ to be the maximizer of~(\ref{s9.e1})
after replacing $\pp_n$ with $\pp_n^{\circ}$. Also let
$\hat{S}_1^{\circ}(\zeta)\equiv\sqrt{n}\pp_n^{\circ}(U_{\zeta,1}^{\tau}
(\hat{\psi}_0^{\circ}),U_{\zeta,2}^{\tau}(\hat{\psi}_0^{\circ})')'$.
Note that the same sample of weights
$\kappa_1^{\circ},\ldots,\kappa_n^{\circ}$ are used for computing
both $\hat{\psi}_0^{\circ}$ and the process
$\{\hat{S}_1^{\circ}(\zeta),\zeta\in[a,b]\}$, so that the proper dependence
between the score statistic and $\hat{\psi}_0$ will be captured.
The structure of the set-up
only requires considering values of $\zeta$ in the set
$\{Y_{(1)},\ldots,Y_{(n)}\}\cap[a,b]$, 
since $\zeta\mapsto\hat{S}_{1}^{\circ}(\zeta)$
does not change over the intervals $[Y_{(j)},Y_{(j+1)})$, $1\leq j\leq n-1$.
Now repeat the bootstrap
procedure a large number of times $\tilde{M}_n$, to obtain
the bootstrapped score processes $\hat{S}_{1,1}^{\circ},\ldots,
\hat{S}_{1,\tilde{M}_n}^{\circ}$. Note that we are allowing
the number of bootstraps to depend on~$n$. Define 
$\zeta\mapsto\hat{\mu}_n(\zeta)\equiv\tilde{M}_n^{-1}\sum_{k=1}^{\tilde{M}_n}
\hat{S}_{1,k}^{\circ}(\zeta)$ and let 
\[\zeta\mapsto\hat{V}_n(\zeta)
=\tilde{M}_n^{-1}\sum_{k=1}^{\tilde{M}_n}\left\{
\hat{S}_{1,k}^{\circ}(\zeta)
-\hat{\mu}_n(\zeta)\right\}\left\{
\hat{S}_{1,k}^{\circ}(\zeta)
-\hat{\mu}_n(\zeta)\right\}'.\]
Now we can compute the test statistics $\hat{T}_n$ and $\tilde{T}_n$
with this choice for $\hat{V}_n$. 

To estimate critical values,
we compute the standardized bootstrap test statistics
$\hat{T}_{n,k}^{\circ}\equiv\sup_{\zeta\in[a,b]}\left\{
\left[\hat{S}_{1,k}^{\circ}(\zeta)-\hat{\mu}_n(\zeta)\right]'
\hat{V}_n^{-1}(\zeta)\left[
\hat{S}_{1,k}^{\circ}(\zeta)-\hat{\mu}_n(\zeta)\right]\right\}$ and
$\tilde{T}_{n,k}^{\circ}\equiv\int_{[a,b]}\left\{
\left[\hat{S}_{1,k}^{\circ}(\zeta)-\hat{\mu}_n(\zeta)\right]'
\hat{V}_n^{-1}(\zeta)\left[
\hat{S}_{1,k}^{\circ}(\zeta)-\hat{\mu}_n(\zeta)\right]\right\}d\zeta$,
for $1\leq k\leq\tilde{M}_n$. For a test of size $\pi$, we compare
the test statistics with the $(1-\pi)$th quantile of the 
corresponding $\tilde{M}_n$ standardized bootstrap statistics.
The reason we subtract off the sample mean when computing
the bootstrapped test statistics is to make sure that we
are approximating the null distribution even when the 
null hypothesis may not be true. What is a little unusual about this
procedure is that the bootstrap must be performed before the 
statistics $\hat{T}_n$ and $\tilde{T}_n$ can be calculated in
the first place. We also reiterate again that we are assuming the 
covariates $Z_i(\cdot)$ are observed at all time points
$V_j\leq V_i$ for which $\delta_j=1$. As noted in section~2, we are
aware that this is not necessarily valid in practice. As pointed out by
a referee this is an important issues and it would be worth investigating
whether the bootstrap weighting scheme could be 
modified to perform and account
for imputation of the missing covariate values. Nevertheless, this issue
is beyond the scope of this paper and we do not pursue it further here.  

\subsection{Asymptotic properties} In this section we establish
the asymptotic validity of the proposed test procedure. 
Let $P$ denote the fixed probability distribution under the null
hypothesis $H_0$, and let $P_n$ be the sequence of probability 
distributions under the contiguous sequence of
alternatives $H_1^n$ defined in C2'.  Note that $P$ and $P_n$
can be equal if $\alpha_{\ast}=0=\eta_{\ast}$. 
We need to study the
proposed procedure under general $P_n$ to determine both its size under
the null and its power under the alternative. 
We will use the notation $\weakpn$ to denote 
weak convergence under $P_n$. We need the following
lemmas and theorem:

\begin{lemma}\label{s9.l1}
The sequence of probability measures $P_n$ satisfies
\begin{eqnarray}
\label{c8.e1}\\
\int \left[ \sqrt{n}(d P_n ^{1/2} - d P ^{1/2})-
\frac{1}{2}\left(U_{\zeta_0,1}^{\tau}(\psi_0^{\ast})(\alpha_{\ast})
+U_{\zeta_0,2}^{\tau}(\psi_0^{\ast})(\eta_{\ast})\right)dP^{1/2}
\right]^2 \rightarrow 0,\nonumber 
\end{eqnarray}
where $\psi_0^{\ast}\equiv(0,0,\beta_0,A_0)$.
\end{lemma}

\begin{lemma}\label{s9.l2}
$\|\hat{\psi}_0-\psi_0^{\ast}\|_{\infty}
\rightarrow 0$ in probability under $P_n$.
\end{lemma}

\begin{theorem}\label{s9.t1}
Under the conditions of section~2, with condition C2 replaced by
C2', $\hat{S}_1$ converges
under $P_n$ in distribution in $l^{\infty}([a,b]^{q+1})$ to the
$(q+1)$-vector process $\zeta\mapsto\mathbb{Z}_{\ast}(\zeta)
+\nu_{\ast}(\zeta)$,
where $\mathbb{Z}_{\ast}$ is a tight, mean zero Gaussian 
$(q+1)$-vector process with
$\mbox{cov}[\mathbb{Z}_{\ast}(\zeta_1),\mathbb{Z}_{\ast}(\zeta_2)]
=\Sigma_{\ast}(\zeta_1,\zeta_2)\equiv
\sigma_{\ast}^{11}(\zeta_1\vee\zeta_2)-\sigma_{\ast}^{12}(\zeta_1)
[\sigma_{\ast}^{22}]^{-1}\sigma_{\ast}^{21}(\zeta_2)$, for all 
$\zeta_1,\zeta_2\in[a,b]$, where, for each $\zeta\in[a,b]$,
\begin{eqnarray*}
\nu_{\ast}(\zeta)&\equiv&\left\{\sigma_{\ast}^{11}
(\zeta\vee\zeta_0)
-\sigma_{\ast}^{12}(\zeta)[\sigma_{\ast}^{22}]^{-1}
\sigma_{\ast}^{21}(\zeta_0)\right\}
\left(\begin{array}{c}\alpha_{\ast}\\ \eta_{\ast}\end{array}\right),\\
\sigma_{\ast}^{11}(\zeta)
&\equiv&\left(\begin{array}{cc}\sigma_{\psi_0^{\ast},\zeta}^{11}
&\sigma_{\psi_0^{\ast},\zeta}^{12}\\ \\ \sigma_{\psi_0^{\ast},\zeta}^{21}
&\sigma_{\psi_0^{\ast},\zeta}^{22}\end{array}\right),\;\;\;\;
\sigma_{\ast}^{12}(\zeta)\;\;\equiv\;\;\left(
\begin{array}{cc}\sigma_{\psi_0^{\ast},\zeta}^{13}
&\sigma_{\psi_0^{\ast},\zeta}^{14}\\ \\ \sigma_{\psi_0^{\ast},\zeta}^{23}
&\sigma_{\psi_0^{\ast},\zeta}^{24}\end{array}\right),\\
\sigma_{\ast}^{21}(\zeta)&\equiv&\left(
\begin{array}{cc}\sigma_{\psi_0^{\ast},\zeta}^{31}
&\sigma_{\psi_0^{\ast},\zeta}^{32}\\ \\ \sigma_{\psi_0^{\ast},\zeta}^{41}
&\sigma_{\psi_0^{\ast},\zeta}^{42}\end{array}\right),\;\;\;\;
\sigma_{\ast}^{22}\;\;\equiv\;\;\left(
\begin{array}{cc}\sigma_{\psi_0^{\ast},\zeta_0}^{33}
&\sigma_{\psi_0^{\ast},\zeta_0}^{34}\\ \\ \sigma_{\psi_0^{\ast},\zeta_0}^{43}
&\sigma_{\psi_0^{\ast},\zeta_0}^{44}\end{array}\right),
\end{eqnarray*}
and where $\sigma_{\theta}^{jk}$, for $1\leq j,k\leq 4$, is as defined
in section~5.2. 
\end{theorem}

The following is the main result on the limiting distribution of the
test statistics. For the remainder of this section, we require 
condition~B4 to hold.  As will be shown in the proof of corollary~\ref{c2},
condition~B4 implies that $V_{\ast}(\zeta)\equiv\Sigma_{\ast}(\zeta,\zeta)$ 
is positive definite
for all $\zeta\in[a,b]$. Note that we will establish consistency of 
$\hat{V}_n$ after we verify the validity of the proposed bootstrap.
\begin{corollary}\label{c2}
Assume~B4 holds and $\hat{V}_n(\zeta)\rightarrow V_{\ast}(\zeta)$
in probability under $P_n$, uniformly over $\zeta\in[a,b]$. 
Then $\hat{T}_n\weakpn\sup_{\zeta\in[a,b]}\left\{
\left[\mathbb{Z}_{\ast}(\zeta)+\nu_{\ast}(\zeta)\right]'\right.$
$\left.\times V_{\ast}^{-1}(\zeta)
\left[\mathbb{Z}_{\ast}(\zeta)+\nu_{\ast}(\zeta)\right]\right\}$
and $\tilde{T}_n\weakpn\int_{[a,b]}\left\{
\left[\mathbb{Z}_{\ast}(\zeta)+\nu_{\ast}(\zeta)\right]'
V_{\ast}^{-1}(\zeta)
\left[\mathbb{Z}_{\ast}(\zeta)\right.\right.$\newline
$\left.\left.+\nu_{\ast}(\zeta)\right]\right\}$.
Thus the limiting null distributions of $\hat{T}_n$ and
$\tilde{T}_n$ are
$\hat{\mathbb{T}}_{\ast}\equiv\sup_{\zeta\in[a,b]}$ $\left\{
\mathbb{Z}_{\ast}'(\zeta)V_{\ast}^{-1}(\zeta)
\mathbb{Z}_{\ast}(\zeta)\right\}$ and 
$\tilde{\mathbb{T}}_{\ast}\equiv\int_{[a,b]}\left\{
\mathbb{Z}_{\ast}'(\zeta)V_{\ast}^{-1}(\zeta)
\mathbb{Z}_{\ast}(\zeta)\right\}d\zeta$, respectively.
\end{corollary}

\begin{remark}\label{r3}
Note that $\nu_{\ast}(\zeta_0)$ equals the matrix 
$\Sigma_{\ast}(\zeta_0,\zeta_0)$ times $(\alpha_{\ast},\eta_{\ast}')'$.
By arguments in the proof of lemma~\ref{l4}, we know that
$\Sigma_{\ast}(\zeta_0,\zeta_0)$ is positive definite. Thus
$\nu_{\ast}(\zeta_0)$ will be strictly nonzero whenever 
$(\alpha_{\ast},\eta_{\ast}')'\neq 0$. Thus both $\hat{T}_n$
and $\tilde{T}_n$ will have power to reject~$H_0$ under
strictly non-null contiguous alternatives $H_1^n$.
\end{remark}

The following theorem is the first step in establishing the validity
of the bootstrap. For brevity, we will use the notation
$\weakpnboot$ to denote conditional convergence of the bootstrap,
either weakly in the sense of corollary~\ref{c1} or in probability, 
but under $P_n$ rather than $P$.
\begin{theorem}\label{s9.t2}
Under the conditions of theorem~\ref{s9.t1}, 
$\hat{S}_1^{\circ}-\hat{S}_1\;\weakpnboot\;\mathbb{Z}_{\ast}$
in $\ell^{\infty}([a,b]^{q+1})$.
\end{theorem}

The following corollary yields the desired consistency of $\hat{V}_n$
and the validity of the proposed bootstrap for obtaining 
critical values. Define 
$\hat{\mathbb{F}}(u)\equiv\tilde{M}_n^{-1}\sum_{k=1}^{\tilde{M}_n}
\ind\left\{\hat{T}_{n,k}^{\circ}\leq u\right\}$ and
$\tilde{\mathbb{F}}(u)\equiv\tilde{M}_n^{-1}\sum_{k=1}^{\tilde{M}_n}
\ind\left\{\tilde{T}_{n,k}^{\circ}\leq u\right\}$.
\begin{corollary}\label{c3}
There exists a sequence $\tilde{M}_n\rightarrow\infty$, as 
$n\rightarrow\infty$, such that 
$\hat{V}_n\weakpn\Sigma_{\ast}$, $\hat{V}_n\weakpnboot\Sigma_{\ast}$, 
and both $\sup_{u\in\re}\left|
\hat{\mathbb{F}}(u)-P\left\{\hat{\mathbb{T}}_{\ast}\leq u\right\}\right|
\weakpnboot 0$ and $\sup_{u\in\re}$ $\left|
\tilde{\mathbb{F}}(u)-P\left\{\tilde{\mathbb{T}}_{\ast}\leq u\right\}\right|
\weakpnboot 0$.
\end{corollary}

\section{Implementation and simulation study} 
We have implemented the proposed estimation
and inference procedures for both the proportional hazards and proportional
odds models. The maximum likelihood estimates were computed using the
profile likelihood $pL_{n}(\zeta)$ defined in section~4. A line search
over the order statistics of $Y$ is used to maximize over $\zeta$, 
while Newton's method is used to maximize over $\psi$. 
The stationary point equation~(\ref{c4:e2}) can be 
used to profile over $A$ for each value of $\zeta$ and $\gamma$. In 
our experience, the computational time of the entire procedure is
reasonable. A thorough simulation study to validate the moderate
sample size performance of this procedure and the proposed bootstrap
procedures of section~8 is underway and will be presented elsewhere.
 
Because of the unusual form of
the statistical tests proposed in section~9, we feel it is worthwhile
at this point to present a small simulation study evaluating their
moderate sample size performance. Both the proportional hazards and
proportional odds models were considered. A single time-independent
covariate with
a standard normal distribution was used, so that $d=q=1$, and the
change-point $Y$ also had a standard normal distribution. The
parameter values were set at $\zeta_0=0$, $\alpha_0=0$, $\beta_0=1$,
$\eta_0\in\{0,-0.5,-1,-2,-3\}$, and $A_0(t)=t$. The range of
$\eta_0$ values includes the null hypothesis $H_0$ (when $\eta_0=0$)
and several alternative hypotheses. The censoring time was
exponentially distributed with rate $0.1$ and truncated at 10.
This resulted in a censoring rate of about 25\%. The sample
size for each simulated data set was 300. For each simulated 
data set, 250 bootstraps were generated with standard exponential
weights truncated at 5, to compute $\hat{V}_n$ and the critical
values for the two test statistics, $\hat{T}_n$ (the ``sup score
test'') and $\tilde{T}_n$ (the ``mean score test''). The range
for $\zeta$ was restricted to the inner 80\% of the $Y$ values.
Each scenario was replicated 250 times.

The results of the simulation study are presented in table~\ref{table1} on
page~\pageref{table1}.
The type~I error (the $\eta_0=0$ column) is quite close to
the targeted 0.05 level, and the power increases with the magnitude
of $\eta_0$. Also, the sup test is notably more powerful than the mean
test for all alternatives. 
We also tried the nonparametric bootstrap and found that
it did not work nearly as well. While it is difficult to make sweeping
generalizations with this small of a numerical study, it appears as if the
proposed test statistics match the theoretical predictions and have
reasonable power. More simulation studies into the properties of these
statistics would be worthwhile, especially studies of the impact of
time-dependent covariates.

\begin{table*}
\caption{Results from the simulation study of the
sup and mean score test statistics in
the proportional hazards and proportional odds models.
The sample size is 300, the level of
censoring approximately 25\%, and the nominal type~I error
is 0.05. 250 replicates were generated for each
configuration. The parameters were set at
$\zeta_0=0$, $\alpha_0=0$, $\beta_0=1$, and $A_0(t)=t$, with
the value of $\eta_0$ varying. The worst-case Monte Carlo standard error
for the power estimates is $0.03=0.50/\sqrt{250}$.}
\label{table1}
\begin{tabular}{c|c|c|c|c|c}
\hline\hline
\multicolumn{6}{c}{Proportional hazards model}\\ \hline
\scriptsize{Sup score test statistic}& Null $\eta_0=0$ & $\eta_0=-0.5$
&  $\eta_0=-1$ & $\eta_0=-2 $ & $\eta_0=-3$ \\ \hline
\scriptsize{mean}&5.078&5.590 &7.874 &13.524&35.507 \\ \hline 
\scriptsize{Standard Deviation}&2.728&2.859 &3.919 &6.992& 11.337 \\\hline
\scriptsize{power}&0.044&0.076 &0.180 &0.536&0.980\\ \hline\hline
\scriptsize{Mean score test statistic}& Null $\eta_0=0$ 
& $\eta_0=-0.5$ & $\eta_0=-1$&  $\eta_0=-2$ &  $\eta_0=-3$ \\\hline
\scriptsize{mean}&1.403 &1.694 &2.560 &5.412 & 5.529 \\ \hline 
\scriptsize{Standard Deviation} &1.206 &1.104 &1.597 &2.492  & 2.683\\ \hline
\scriptsize{power}&0.040 &0.050 &0.120 &0.236  &0.304 \\ \hline \hline
\multicolumn{6}{c}{Proportional odds model}\\ \hline
\scriptsize{Sup score test statistic}& Null $\eta_0=0$ & $\eta_0=-0.5$
&  $\eta_0=-1$ & $\eta_0=-2 $ & $\eta_0=-3$ \\ \hline
\scriptsize{mean}
&3.950&4.762 &5.693 & 8.327  &13.956 \\\hline \scriptsize{Standard
Deviation}&2.390&1.610 & 1.255 &2.901 & 4.244\\ \hline
\scriptsize{power}&0.043&0.068 &0.112 &0.364 &0.660 \\ \hline\hline
\scriptsize{Mean score test statistic}& Null 
$\eta_0=0$ & $\eta_0=-0.5$ & $\eta_0=-1$&  $\eta_0=-2$ 
&  $\eta_0=-3$ \\ \hline
\scriptsize{mean}&1.177 &1.912 &2.848 &3.265 &4.349 \\ \hline 
\scriptsize{Standard Deviation} &0.946 & 1.078 &1.360 &1.498  &1.718 \\ \hline
\scriptsize{power}&0.048 &0.056 &0.116 &0.167  &0.285 \\ \hline\hline
\end{tabular}
\end{table*}

\section{Proofs} 
       
{\it Proof of lemma~\ref{l.v1}.} Verification of~D1 is straightforward.
For~D2, we have for all $u\geq 0$,
\[\left|\frac{\ddot{\Lambda}(u)}{\dot{\Lambda}(u)}\right|=
\frac{\Exp{W^2e^{-uW}}}{\Exp{We^{-uW}}}\leq\frac{\Exp{W^2}}{\Exp{W}}<\infty.\]
The second-to-last inequality requires some justification. Note that the
probability measure $Qf(W)\equiv\Exp{f(W)W}/\Exp{W}$ is well-defined 
for functions $f$ bounded by $O(W^3)$ by the positivity of $W$ and the 
existence of a fourth moment. Now we have
\[\frac{\Exp{W^2e^{-uW}}}{\Exp{We^{-uW}}}=\frac{Q[We^{-uW}]}
{Q[e^{-uW}]}\leq Q[W]=\frac{\Exp{W^2}}{\Exp{W}},\]
since $e^{-uW}$ uniformly down-weights larger values of $W$ and thus forces 
the left term of the inequality to be decreasing in~$u$. This proves
the first part.

For the second part, take $c_0=c$, and note that
\[|u^c\Lambda(u)|=\Exp{u^ce^{-uW}}=\Exp{W^{-c}(uW)^ce^{-uW}}
\leq k\Exp{W^{-c}},\]
where $k=\sup_{x\geq 0}x^ce^{-x}=c^ce^{-c}<\infty$. Similarly,
\[|u^{1+c}\dot{\Lambda}(u)|=\Exp{u^{1+c}We^{-uW}}
=\Exp{W^{-c}(uW)^{1+c}e^{-uW}}\leq k'\Exp{W^{-c}},\]
where $k'=\sup_{x\geq 0}x^{1+c}e^{-x}=(1+c)^{1+c}e^{-1-c}<\infty$.
This concludes the proof.$\Box$

{\it Proof of lemma~\ref{l1}.} Suppose that
\begin{eqnarray}
\;\;\;\;G\left(\int_0^t \tilde{Y}(s)e^{r_{\xi}(u;Z,Y)}dA(u) \right) = G
  \left(\int_0^t \tilde{Y}(s)e^{r_{\xi_0}(u;Z,Y)}dA_0(u)
  \right)~\label{c12:e1}
\end{eqnarray}
for all $t \in [0, \tau]$ almost surely under $P$. The target is
to show that~(\ref{c12:e1}) implies that $\xi = \xi_0$ and $A =
A_0$ on $[0, \tau]$. By condition~A1, (\ref{c12:e1})~implies
\[\int_0^t e^{r_{\xi}(u;Z,Y)}dA(u)=
\int_0^t e^{r_{\xi_0}(u:Z,Y)}dA_0(u)\]
for all $t \in [0, \tau]$ almost surely. 
Taking the Radon-Nikodym derivative of both
sides with respect to $A_0$, and taking logarithms, we obtain
\begin{eqnarray}
&&\beta'Z(t)+(\alpha+\eta'Z_2(t))\ind\{Y>\zeta\}-\beta_0'Z(t)\label{c12:e2}\\
&&\mbox{\hspace{1.5in}}
-(\alpha_0+\eta_0'Z_2(t))\ind\{Y>\zeta_0\} +\log(\tilde{a}(t))=0,\nonumber
\end{eqnarray}
almost surely, where $\tilde{a} \equiv dA/dA_0$. 

Assume that $\zeta>\zeta_0$. Now choose $y<\zeta_0$ such
that $y\in\tilde{V}(\zeta_0)$ and
$\mbox{var}[Z(t_1)|Y=y]$ is positive definite, where $t_1$ is
as defined in~B3.  Note that this is possible by assumptions~B2
and~B3. Conditioning the left-hand side of~(\ref{c12:e2}) 
on $Y=y$ and evaluating at $t=t_1$ yields that $\beta=\beta_0$.
Now choose $\zeta_0<y<\zeta$ such that $y\in\tilde{V}(\zeta_0)$ and
$\mbox{var}[Z(t_2)|Y=y]$ is positive definite. Conditioning
the left-hand side of~(\ref{c12:e2}) on $Y=y$, and evaluating
at $t=t_2$ yields that $\eta_0=0$. Because the density of $Y$ is
positive in $\tilde{V}(\zeta_0)$, we also see that $\alpha_0=0$.
But this is not possible by condition~C2. A similar argument can 
be used to show that $\zeta<\zeta_0$ is impossible. Thus $\zeta=\zeta_0$.
Now it is not hard to argue that condition~B3 forces
$\beta=\beta_0$, $\eta=\eta_0$ and $\alpha=\alpha_0$.
Hence $\log(\tilde{a}(t))=0$ for all $t\in[0,\tau]$, 
and the proof is complete.$\Box$

{\it Proof of lemma~\ref{l2}.} Note that for each $n$, maximizing
the log-likelihood
over $A$ is equivalent to maximizing over a fixed number of parameters
since the number of jumps $K\leq n$.  Thus maximizing over the
whole parameter $\theta$ involves maximizing an empirical average of
functions that are smooth over $\psi$ and cadlag over $\zeta$. 
Note also that 
\[\|\hat{A}_n-A_0\|_{[0,\tau]}=\sum_{j=1}^{K}
\left(\left|\hat{A}_n(T_j-)-A_0(T_j)\right|\vee
\left|\hat{A}_n(T_j)-A_0(T_j)\right|\right),\]
where $\|\cdot\|_B$ is the uniform norm over the set $B$, 
and thus $\|\hat{A}_n-A_0\|_{[0,\tau]}$ is measurable.
Hence the uniform distance between 
$\hat{\theta}_n$ and $\theta_0$ is also measurable. 
Thus almost sure convergence of
$\hat{\theta}_n$ is equivalent to outer almost sure convergence.
Now we return to the proof. Assume
\begin{eqnarray}
\lim\sup_{n \rightarrow \infty} \hat{A}_n(\tau) =
\infty,\label{c12:e4}
\end{eqnarray}
with probability $>0$. 
We will show that this leads to a contradiction. It is now possible to 
choose a data sequence such that~(\ref{c12:e4}) holds and
$\tilde{G}_n\equiv\pp_n N\rightarrow \tilde{G}_0\equiv P_0 N$ 
uniformly, since the latter happens with probability~1. Fix one such
sequence $\{n\}$, and define $\theta_n=(\xi_0,A_n)$, where $A_n=\tilde{G}_n$. 
Note that the log-likelihood difference, 
$\tilde{L}_n(\hat\theta_n)-\tilde{L}_n(\theta_n)$,
should be non-negative for all $n$, since $\hat{\theta}_n$
maximizes the log-likelihood. We are going to show that the
difference is asymptotically negative under the
assumption~(\ref{c12:e4}).

Now choose a subsequence $\{n_k\}$ such that $\hat{A}_{n_k}(\tau)
\rightarrow \infty$, as $k\rightarrow\infty$. We now have,
for $c_0>0$ from assumption~D2, that
$L_{n_k}(\hat{\theta}_{n_k})-L_{n_k}(\theta_{n_k})$
\begin{eqnarray}
\mbox{\hspace{0.7in}} &\le& O(1) 
+\mathbb{P}_{n_k} \delta \left[\log 
\left(n_k \Delta \hat{A}_{n_k}(V) \right)+\log\left(
-\dot{\Lambda}(H^{\hat{\theta}_n}(V))\right)\right]\nonumber\\
&&-\mathbb{P}_{n_k}(1-\delta)G(H^{\hat{\theta}_n}(V))\nonumber\\
&\leq& O(1)+\mathbb{P}_{n_k} \delta\log 
\left(n_k \Delta \hat{A}_{n_k}(V) \right)
-\mathbb{P}_{n_k}(\delta+c_0)\log\hat{A}_n(V),\label{c12:e5}
\end{eqnarray}
since, for all $u>0$, 
$\log\dot{G}(u)=\log[-\dot{\Lambda}(u)]-\log[\Lambda(u)]$;
$\log[-\dot{\Lambda}(u)]=\log[-u^{1+c_0}$
$\dot{\Lambda}(u)]-(1+c_0)\log(u)
\leq O(1)-(1+c_0)\log(u)$ by condition~D2; and since
$\log\Lambda(u)=\log[u^{c_0}\Lambda(u)]-c_0\log(u)\leq O(1)-c_0\log(u)$
also by condition~D2.

Next we take a partition of $[0, \tau]$, $0=v_0<v_1<\cdots<v_M=\tau$,
for some finite $M$. The right hand side of~(\ref{c12:e5}) is now
dominated by
\begin{eqnarray}\label{c12:e6}
O(1)+\log \hat{A}_{n_k}(\tau) \mathbb{P}_{n_k} \left(\delta 
\ind\{V \in [v_{M-1}, \infty] \} - (\delta + c_0)\ind\{V \in
[\tau, \infty \}  \right)&&\\
+\sum_{m=1}^{M-1}\log \hat{A}_{n_k}(v_m) \mathbb{P}_{n_k} 
\left(\delta \ind\{V \in [v_{m-1}, v_m] \} - (\delta + c_0)
\ind\{V \in [v_{m}, v_{m+1}] \}  \right).&&\nonumber
\end{eqnarray}
For a fixed constant $c > 1$, we can choose this partition such that
\begin{eqnarray*}
P_0N(\tau)\ind\{V \in [v_{M-1}, \infty]\} = P_0[N(\tau) +
c_0/c]\ind\{V \in [\tau, \infty] \},
\end{eqnarray*}
and, for $m = 1,...,M-1,$
\begin{eqnarray*}
P_0N(\tau)\ind\{V \in [v_{m-1}, v_m]\} = P_0[N(\tau) +
c_0/c]\ind\{V \in [v_m, v_{m+1}] \}.
\end{eqnarray*}
Recalling that $\tilde{G}_n\rightarrow \tilde{G}_0$ uniformly, we
obtain that~(\ref{c12:e6}) tends to $-\infty$ as $k \rightarrow \infty$,
which is the intended contradiction. Thus,  $\limsup_{n
\rightarrow \infty} \hat{A}_n(\tau) < \infty$ almost surely.$\Box$

{\it Proof of theorem~\ref{t1}.} 
By the opening arguments in the proof of lemma~\ref{l2}, we have that
outer almost sure convergence is equivalent to the usual almost
sure convergence in this instance. Note that $\{\hat{A}_n(\tau)\}$ is bounded
almost surely, $\tilde{G}_n\rightarrow\tilde{G}_0$ almost
surely, and the class
\[{\cal F}_{(k)}\equiv\left\{W(t;\theta):t\in[0,\tau],\xi\in{\cal X},
A\in{\cal A}_{(k)}\right\},\]
where ${\cal A}_{(k)}\equiv\{A\in{\cal A}:A(\tau)\leq k\}$,
is Donsker (and hence also Glivenko-Cantelli) for every $k<\infty$
by lemma~\ref{l.t1.1} below. By similar arguments to those used
in lemma~\ref{l.t1.1}, we have that the class
$\{G(H^{\theta}(V)):\xi\in{\cal X},A\in{\cal A}_{(k)}\}$ is also
Glivenko-Cantelli for all $k<\infty$. We therefore
have the following with probability~1:
$\{\hat{A}_n(\tau)\}$ is bounded asymptotically, $\tilde{G}_n\rightarrow
\tilde{G}_0$ uniformly, $(\mathbb{P}_n-P)W(\cdot;\hat{\theta}_n)
\rightarrow 0$ uniformly, and $(\mathbb{P}_n-P)\left[
G(H^{\hat{\theta}_n}(V))-G(H^{\theta_n}(V))\right]\rightarrow 0$. 
Now, fix a sequence $\{n\}$ for which these last four asymptotic events
hold. We can now use the Helly selection theorem to find a
subsequence $\{n_k\}$ and a function $A$ such that
$\hat{A}_{n_k}(t) \rightarrow A(t)$ for all $t \in [0, \tau]$ at
which $A$ is continuous. From~(\ref{c4:e2}), we obtain
\[|\hat{A}_{n_k}(s)-\hat{A}_{n_k}(t)| \le
O(1)\mathbb{P}_{n_k}|N(s)-N(t)|\rightarrow O(1)|\tilde{G}_0(s)
-\tilde{G}_0(t)|,\] 
for all $s,t \in [0, \tau]$. Since $\tilde{G}_0$ is continuous by
condition~C3, we know that $A$ must be continuous on all of $[0,\tau]$.
Thus $\hat{A}_{n_k} \rightarrow A$ uniformly. 
Without loss of generality, we can also assume that along this
subsequence $\hat{\xi}_{n_k} \rightarrow \xi$ for some $\xi \in
{\cal X}\equiv\Upsilon\times B_1\times B_2\times(a,b)$. Denote 
$\theta=(\xi,A)$.

Consider now $\theta_n \equiv (\xi_0, A_n)$, where
\begin{eqnarray*}
A_n(t) \equiv \int_0^t \frac{d\tilde{G}_n(u)}{PW(u;
\theta_0)}.
\end{eqnarray*}
We can use the same technique as in the derivation
of~(\ref{c4:e2}) to show that $A_0$ satisfies
\begin{eqnarray*}
A_0(t)\equiv\int_0^t \frac{d\tilde{G}_0(u)}{PW(u;\theta_0)},
\end{eqnarray*}
for all $t\in[0,\tau]$. Thus $A_{n_k}\rightarrow A_0$ uniformly,
as $k\rightarrow\infty$. At this point, we have
\begin{eqnarray*}
0&\leq&\tilde{L}_{n_k}(\hat{\theta}_{n_k})-\tilde{L}_{n_k}(\theta_{n_k})\\
&=&\int_0^{\tau}\log\left[\frac{PW(u;\theta_0)}
{\mathbb{P}_{n_k}W(u;\hat{\theta}_{n_k})}\right]
d\tilde{G}_{n_k}(u)-\mathbb{P}_{n_k}\left[G(H^{\hat{\theta}_{n_k}}(V))-
G(H^{\theta_{n_k}}(V))\right]\\
&\rightarrow&\int_0^{\tau}\log\frac{dA(u)}{dA_0(u)}d\tilde{G}_0(u)
-P\left[G(H^{\theta}(V))-G(H^{\theta_0}(V))\right]\\
&=&\int\log\frac{dP_{\theta}}{dP}dP\\
&\leq&0.
\end{eqnarray*}
But this forces $\theta=\theta_0$ by the identifiability of
the model as given in lemma~\ref{l1}. Thus all convergent 
subsequences of $\hat{\theta}_n$,
on a set of probability~1, converge to $\theta_0$. The desired
result now follows.$\Box$

\begin{lemma}\label{l.t1.1}
$\forall k<\infty$, the class
${\cal F}_{(k)}\equiv\left\{W(t;\theta):t\in[0,\tau],\xi\in{\cal X},
\right.$ $\left.A\in{\cal A}_{(k)}\right\}$,
is $P$-Donsker.
\end{lemma}

{\it Proof.} Routine arguments can be used to establish that
the class ${\cal F}_1\equiv
\{e^{r_{\xi}(t;Z,Y)}:t\in[0,\tau],\xi\in{\cal X}\}$ is
Donsker. Consider the map 
\[h\in D[0,\tau]\mapsto
\left\{\int_0^th(s)dA(s):t\in[0,\tau],A\in{\cal A}_{(k)}\right\}
\in\ell^{\infty}([0,\tau]\times{\cal A}_{(k)}),\]
and note that it is uniformly equicontinuous and linear.
Thus the class 
\[{\cal F}_2\equiv\left\{\int_0^t e^{r_{\xi}(s;Z,Y)}dA(s):
t\in[0,\tau],\xi\in{\cal X},A\in{\cal A}_{(k)}\right\}\]
is Donsker by the continuous mapping theorem.
Now condition~D1 ensures that both $\dot{G}$ and
$\ddot{G}/\dot{G}$ are Lipschitz on compacts.  This fact,
combined with the facts that sums of Donsker classes 
are Donsker and products of bounded Donsker classes are
Donsker, yields the desired results.$\Box$

{\it Proof of lemma~\ref{l3}.} By the smoothness assumed in~D1 of the involved
derivatives, we have for each $\zeta\in[a,b]$ and $\psi^{\ast}\in\Psi$, 
\[\lim_{t\downarrow 0}\sup_{h^{\ast}\in\mbox{lin}\,\Psi:
\rho_1(h^{\ast})\leq 1}
\sup_{h\in{\cal H}_r}\left|\int_0^1 h^{\ast}
\left(\sigma_{\psi^{\ast}+st h^{\ast}}(h)-\sigma_{\psi^{\ast}}(h)\right)ds
\right|=0.\]
Thus, $\sup_{h\in{\cal H}_r}\left|PU_{\zeta}^{\tau}(\psi^{\ast}+h^{\ast})(h)
-PU_{\zeta}^{\tau}(\psi^{\ast})(h)+h^{\ast}\left(\sigma_{\psi^{\ast}}(h)
\right)\right|=o(\rho_1(h^{\ast}))$, as $\rho_1(h^{\ast})\rightarrow 0$.$\Box$

{\it Proof of lemma~\ref{l4}.}
First note that for any $h=(h_1,h_2,h_3,h_4)\in{\cal H}_{\infty}$,
$\sigma_{\theta_0}(h)=\mathbb{A}(h)+\mathbb{B}(h)$, where
$\mathbb{A}(h)=\left(h_1,h_2,h_3,g_0h_4\right)$, 
$\mathbb{B}(h)=\sigma_{\theta_n}(h)-\mathbb{A}(h)$,
and $g_0(u)=P\left[\tilde{Y}(u)e^{r_{\xi_0}(u;Z,Y)}
\hat{\Xi}_{\theta_0}^{(0)}(\tau)\right]$. It is not hard to verify that
since $g_0$ is bounded below, $\mathbb{A}$ 
is one-to-one and onto with continuous 
inverse defined by $\mathbb{A}^{-1}(h)=(h_1,h_2,h_3,h_4/g_0)$.  
It is also not hard to
verify that the operator $\mathbb{B}$ 
is compact as an operator on ${\cal H}_r$
for any $0<r<\infty$. Thus the first part of the theorem is proved 
by lemma~25.93 of \cite{v98}, if
we can show that $\sigma_{\theta_0}$ is one-to-one. This will then imply
that for each $r>0$, there is an $s>0$ with 
$\sigma_{\theta_0}^{-1}({\cal H}_s)\subset{\cal H}_r$. Now we have
\[\inf_{\psi\in\mbox{lin}\,\Psi}
\frac{\|\psi(\sigma_{\theta_0}(\cdot))\|_{(r)}}
{\|\psi\|_{(r)}}\geq \inf_{\psi\in\mbox{lin}\,\Psi}
\frac{\sup_{h\in\sigma_{\theta_0}^{-1}({\cal H}_s)}
|\psi(\sigma_{\theta_0}(h))|}{\|\psi\|_{(r)}}=
\inf_{\psi\in\mbox{lin}\,\Psi}\frac{\|\psi\|_{(s)}}{\|\psi\|_{(r)}}\]
$\geq s/(4r)$, since $\|\psi\|_{(r)}\leq 4(r/s)\|\psi\|_{(s)}$.
Thus $\psi\mapsto\psi(\sigma_{0}(\cdot))$ is continuously invertible
on its range by proposition~A.1.7 of \cite{bkrw98}. That it is also onto with
inverse $\psi\mapsto\psi(\sigma_{\theta_0}^{-1})$ follows from
$\sigma_{\theta_0}$ being onto. All that remains is verifying that
$\sigma_{\theta_0}$ is one-to-one.

Let $h \in \mathcal{H}_{\infty}$ such that $\sigma_{\theta_0}(h)=0$. 
For the one-dimensional submodel defined by the map $s \rightarrow\psi_{0s} 
\equiv \psi_0 + s(h_1, h_2, h_3, \int_0^{(\cdot)}h_4(u)dA_0(u))$, we have
\begin{eqnarray}
P \{ \frac{ \partial}{\partial s}L_1(\psi_{0s},\zeta_0)|_{s=0}\}^2 = P
\{U^{\tau}_{\zeta_0}(\psi_0)(h)\}^2=0.~\label{c12:e9}
\end{eqnarray}
Define the random set
$\mathcal{S}(n,\tilde{y},t) \equiv \{(N,\tilde{Y}): N(u) = n(u), 
\tilde{Y}(u) = \tilde{y}(u), u \in [t, \tau] \}$.
The equality~(\ref{c12:e9}) implies that
$P\{U_{\zeta_0}^{\tau}(\psi_0)(h)|\mathcal{S}(n,y,t)\}^2=0$
for all $\mathcal{S}$ such that $P\{\mathcal{S}(n,y,t)\} > 0$, which
implies that $U^t_{\zeta_0}(\psi_0)(h)=0$ almost surely for all $t \in [0,
\tau]$. Consider the set on which the observation $(X, \delta, Z,
Y)$ is censored at a time $t \in [0, \tau]$. From (\ref{c12:e9})
and the preceding argument,
\begin{eqnarray}
R_{\zeta_0,\psi_0}^t(h_1\ind(Y > \zeta_0) +
h_2'Z_2(t)\ind(Y > \zeta_0)+h_3'Z(t)+h_4)=0.~\label{c12:e11}
\end{eqnarray}
Taking the Radon-Nikodym derivative of (\ref{c12:e11}) with respect to $A_0$
and dividing throughout by $e^{r_{\xi_0}(t;Z,Y)}$ yields
\begin{eqnarray}
\tilde{Y}(t)(h_1\ind(Y > \zeta_0) +
h_2'Z_2(t)\ind(Y > \zeta_0)+h_3'Z(t)+h_4(t))=0.~\label{c12:e12}
\end{eqnarray}
Arguments quite similar to those used in the proof of lemma~\ref{l1}
can now be used to verify that~(\ref{c12:e12}) forces $h=0$.
Hence $\sigma_{\theta_0}(h)=0$ implies $h=0$, and thus
$\sigma_{\theta_0}$ is one-to-one.$\Box$

{\it Proof of lemma~\ref{l5}.} For the first part, note that 
$t\mapsto\tilde{Y}(t)$ has total variation bounded by~1; 
and, by the model assumptions, the total variation of 
$t\mapsto e^{r_{\xi}(t;Z,Y)}$ is bounded by a universal constant that doesn't
depend on $\theta$. Thus there exists a universal constant $k_{\ast}$
such that $\|\mathbb{P}_n W(\cdot;\hat{\theta}_n)\|_v\leq
k_{\ast}\mathbb{P}_n|\hat{\Xi}^{(0)}_{\hat{\theta}_n}|$.  By the smoothness of
the functions involved, and the fact that $u\mapsto\log(u)$ is Lipschitz
on compacts bounded above zero, we obtain the first result of the lemma.
The consistency part follows from 
lemma~\ref{l.t1.1} combined with theorem~\ref{t1}, the 
continuity of $\theta\mapsto PW(\cdot;\theta)$, and reapplication
of the Lipschitz continuity of $u\mapsto\log(u)$.$\Box$

{\it Proof of lemma~\ref{l6}.} The right-hand derivative of 
$P(L_1(\psi,\zeta))$ with respect to $\zeta$ at $\zeta=\zeta_0$ is:
$\left.(\partial^{+}/(\partial\zeta)) 
P(L_1(\psi, \zeta))\right|_{\zeta=\zeta_0}$
\begin{eqnarray*}
&=&\int\left 
\{P[l_1^{\psi}(V,\delta,Z)|Y=y+]-P[l_2^{\psi}(V,\delta,Z)|Y=y+] \right \} 
\tilde{\delta}_{\zeta_0}(y)\tilde{h}(y)dy \\
&=&\left(P[l_1^{\psi}(V,\delta,Z) |Y=\zeta_0+]-
P[l_2^{\psi}(V,\delta,Z)|Y=\zeta_0+]\right)\tilde{h}(\zeta_0),
\end{eqnarray*}
where the superscript~$+$ denotes differentiating from the right
and $\tilde{\delta}_{\zeta_0}(y)$ is the Dirac delta function
assigning counting measure~1 to the event $\{y=\zeta_0\}$. Now,
$P[l_1^{\psi}(V,\delta,Z)|Y=\zeta_0+]-
P[l_2^{\psi}(V,\delta,Z)|Y=\zeta_0+]$
\[=\int\left[ l_1^{\psi}(v,d,z)-l_2^{\psi}(v,d,z)\right]
\ell_2(v,d,z)\ell_0^{+}(v,d,z)d\mu(v,d,z)\]
$\equiv\tilde{R}^{+}(\psi)$,
where $\ell_j(v,d,z)\equiv\exp\{l_j^{\psi_0}(v,d,z)\}$, for $j=1,2$;
$\mu(v,d,z)$ is the dominating measure; and $\ell_0^{+}(v,d,z)$ 
consists of the remaining components of the conditional distribution
of $(V,\delta,Z)$ given $Y=\zeta_0+$. Note that under the model
assumptions, $\ell_0^{+}$ does not depend on the parameters. Thus
\begin{eqnarray*}
\tilde{R}^{+}(\psi_0)&=&\int\left[ l_1^{\psi_0}(v,d,z) 
-l_2^{\psi_0}(v,d,z)\right]\ell_2(v,d,z)\ell_0^{+}(v,d,z)d\mu(v,d,z)\\
&=&\int\log\left[\frac{\ell_1\ell_0^{+}}{\ell_2\ell_0^{+}}\right]
\ell_2\ell_0^{+}d\mu
\;\;<\;\;\log\int\left[\frac{\ell_1\ell_0^{+}}
{\ell_2\ell_0^{+}}\right]\ell_2\ell_0^{+}d\mu\\
&=&\log\int\ell_1(v,d,z)\ell_0^{+}(v,d,z)d\mu(v,d,z)\;\;=\;\;0,
\end{eqnarray*}
since the integral of a density is~1. Thus 
$\dot{X}_{\zeta_0}^{+}(\gamma_0,\Gamma_0)<0$. 

A similar argument is
used for the left-hand derivative. In this case, the true density
of $(V,\delta,Z)$ given $Y=\zeta_0$ is $\ell_1^{\psi_0}(v,d,z)
\ell_0^{-}(v,d,z)$, where $\ell_0^{-}$ does not involve the parameters.
We now have 
\begin{eqnarray*}
\lefteqn{P[l_1^{\psi}(V,\delta,Z)|Y=\zeta_0]-
P[l_2^{\psi}(V,\delta,Z)|Y=\zeta_0]}\mbox{\hspace{1.0cm}}&&\\
&=&\int\left[ l_1^{\psi_0}(v,d,z) 
-l_2^{\psi_0}(v,d,z)\right]\ell_2(v,d,z)\ell_0^{-}(v,d,z)d\mu(v,d,z)\\
&=&-\int\log\left[\frac{\ell_2\ell_0^{-}}{\ell_1\ell_0^{-}}\right]
\ell_1\ell_0^{-}d\mu
\;\;>\;\;-\log\int\left[\frac{\ell_2\ell_0^{-}}{\ell_1\ell_0^{-}}\right]
\ell_1\ell_0^{-}d\mu\\
&=&\log\int\ell_2(v,d,z)\ell_0^{-}(v,d,z)d\mu(v,d,z)\;\;=\;\;0,
\end{eqnarray*}
and thus we conclude that $\dot{X}_{\zeta_0}^{-}(\gamma_0,\Gamma_0)>0$.$\Box$

{\it Proof of lemma~\ref{l7}.} This follows from lemma~\ref{l6}, the local
concavity of $\tilde{X}$, and the 
smoothness of the derivatives involved.$\Box$

{\it Proof of lemma~\ref{l8}.} 
Note that $\tilde{X}_n(\zeta,\eta,\Gamma)$
\[=\mathbb{P}_n\left[-\int_0^{\tau}\left\{\Gamma(t)-\Gamma_0(t)\right\}dN(t)
+\tilde{W}(\zeta,\eta,A_n^{(\Gamma)})
-\tilde{W}(\zeta_0,\eta_0,A_n^{(\Gamma_0)})\right],\] 
where 
$\tilde{W}(\zeta,\gamma,A)\equiv l_1^{\psi}(V,\delta,Z)\ind\{Y\leq\zeta\}
+l_2^{\psi}(V,\delta,Z)\ind\{Y>\zeta\}$. The classes
\[\left\{\int_0^{\tau}\left\{\Gamma(t)-\Gamma_0(t)\right\}dN(t):
\|\Gamma-\Gamma_0\|_{\infty}\leq\epsilon,\|\Gamma\|_v\leq k_0\right\},\] 
for any $\epsilon>0$, and
$\left\{\tilde{W}(\zeta,\lambda):(\zeta,\lambda)
\in B_{\epsilon_2}^{\ast k_0}\right\}$, for some $\epsilon_2>0$,
can be shown to be Donsker. That this holds for the second class
follows from arguments similar to those used in the proof
of lemma~\ref{l.t1.1}. For the first class, note that
$\int_0^{\tau}\Gamma(t)dN(t)=\delta\Gamma(V)$. Since $\|\Gamma\|_v\leq k_0$,
$\Gamma$ can be written as the 
difference between two monotone increasing functions,
each with total variation bounded by $k_0$. By theorem~2.7.5 of \cite{vw96}, 
the class of all monotone functions with a given compact range is universally
Donsker. Since sums of Donsker classes are Donsker, we have that the class
$\{\Gamma(V):\|\Gamma\|_v\leq k_0\}$ is Donsker. That the first class
is Donsker now follows since products of bounded Donsker classes are Donsker.
Since we also have that
$\sqrt{n}(\tilde{G}_n-\tilde{G}_0)$ converges to a Gaussian process,
we have that 
\[\sqrt{n}(\mathbb{P}_n-P)
\left[-\int_0^{\tau}\left\{\Gamma(t)-\Gamma_0(t)\right\}dN(t)
+\tilde{W}(\zeta,\eta,A_n^{(\Gamma)})
-\tilde{W}(\zeta_0,\eta_0,A_n^{(\Gamma_0)})\right]\]
converges weakly in $\ell^{\infty}(B_{\epsilon_2}^{\ast k_0})$
to the tight Gaussian process 
\[\mathbb{G}\left[-\int_0^{\tau}\left\{\Gamma(t)-\Gamma_0(t)\right\}dN(t)
+\tilde{W}(\zeta,\eta,A_0^{(\Gamma)})
-\tilde{W}(\zeta_0,\eta_0,A_0^{(\Gamma_0)})\right],\]
where $\mathbb{G}$ is the Brownian bridge measure.

By the smoothness of the functions and derivatives involved, 
we also have
$\sqrt{n}\left\{P\left[-\int_0^{\tau}\left\{\Gamma(t)-\Gamma_0(t)\right\}dN(t)
+\tilde{W}(\zeta,\eta,A_n^{(\Gamma)})
-\tilde{W}(\zeta_0,\eta_0,A_n^{(\Gamma_0)})\right]-\right.$
\newline $\left.\tilde{X}(\zeta,\eta,\Gamma)\right\}\;\;=\;\;
\sqrt{n}P\left[\tilde{W}(\zeta,\eta,A_n^{(\Gamma)})
-\tilde{W}(\zeta_0,\eta_0,A_n^{(\Gamma_0)})
-\tilde{W}(\zeta,\eta,A_0^{(\Gamma)})\right.$ \newline
$\left.+\tilde{W}(\zeta_0,\eta_0,A_0)\right]\;\;=\;\;$
$-\sqrt{n}\int_0^{\tau}\left\{P[W(t;\theta_0(\zeta,\lambda))]e^{-\Gamma(t)}
-P[W(t;\theta_0)]e^{-\Gamma_0(t)}\right\}$ $\times
\left[d\tilde{G}_n(t)-d\tilde{G}_0(t)\right]+\epsilon_n(\zeta,\lambda)
\equiv-\int_0^{\tau}\tilde{C}(t;\zeta,\lambda)d{\cal Z}_n(t)+
\epsilon_n(\zeta,\lambda)$,
where $\|\epsilon_n\|_{\infty}$ $=o_P(1)$. The fact that the class
of functions $\{\tilde{C}(\cdot;\zeta,\lambda):(\zeta,\lambda)\in
B_{\epsilon_2}^{\ast k_0}\}$ has uniformly bounded total variation yields 
asymptotic linearity and normality of $\left\{
\int_0^{\tau}\tilde{C}(t;\zeta,\lambda)d{\cal Z}_n(t):(\zeta,\lambda)
\in B_{\epsilon_2}^{\ast k_0}\right\}$,
and the desired result follows.$\Box$

{\it Proof of theorem~\ref{t.l9}.}  By lemma~\ref{l8}, 
\[-\tilde{X}(\hat{\zeta}_n,\hat{\gamma}_n,\hat{\Gamma}_n)
=(\tilde{X}_n-\tilde{X})(\hat{\zeta}_n,\hat{\gamma}_n,\hat{\Gamma}_n)
-\tilde{X}_n(\hat{\zeta}_n,\hat{\gamma}_n,\hat{\Gamma}_n)\leq O_P(n^{-1/2}).\]
Combining this with lemma~\ref{l7}, we obtain
$\sqrt{n}|\hat{\zeta}_n-\zeta_0|$
\begin{eqnarray*}
&=&\sqrt{n}|\hat{\zeta}_n-\zeta_0|\ind\{(\hat{\zeta}_n,\hat{\gamma}_n,
\hat{\Gamma}_n)\in B_{\epsilon_1}^{\ast k_0}\}+
\sqrt{n}|\hat{\zeta}_n-\zeta_0|\ind\{(\hat{\zeta}_n,\hat{\gamma}_n,
\hat{\Gamma}_n)\not\in B_{\epsilon_1}^{\ast k_0}\}\\
&\leq&-\sqrt{n}k_1^{-1}\tilde{X}(\hat{\zeta}_n,\hat{\gamma}_n,
\hat{\Gamma}_n)+o_P(1)\\
&\leq& O_P(1).
\end{eqnarray*}
Thus the first part of the lemma is proved.

For the second part, denote $U_{0\zeta}^{\tau}(\psi)\equiv
P U_{\zeta}^{\tau}(\psi)$. By arguments similar to those used
in the proof of lemma~\ref{l.t1.1}, we can verify that for some $e_1>0$, 
${\cal F}\equiv
\{U_{\zeta}^{\tau}(\psi)(h):\|\theta-\theta_0\|_{\infty}\leq e_1,
h\in{\cal H}_1\}$ is Donsker. Moreover, the continuity of the
functions involved also yields that, as 
$\|\theta-\theta_0\|_{\infty}\rightarrow 0$,
$\sup_{h\in{\cal H}_1}
P\left(U_{\zeta}^{\tau}(\psi)(h)-U_{\zeta_0}^{\tau}(\psi_0)(h)\right)^2
\rightarrow 0$. Thus
\begin{eqnarray}
\sqrt{n}\left(U_{n\hat{\zeta}_n}^{\tau}(\hat{\psi}_n)
-U_{0\hat{\zeta}_n}^{\tau}(\hat{\psi}_n)-U_{n\zeta_0}^{\tau}(\psi_0)
+U_{0\zeta_0}^{\tau}(\psi_0)\right)&=&o_P^{{\cal H}_1}(1).\label{l9.e1}
\end{eqnarray}
Note also that $\sqrt{n}|\hat{\zeta}_n-\zeta_0|=O_P(1)$ implies that
$\sqrt{n}\left(U_{0\hat{\zeta}_n}^{\tau}(\hat{\psi}_n)
-U_{0\zeta_0}^{\tau}(\hat{\psi}_n)\right)=o_P^{{\cal H}_1}(1)$.
Thus, since
$U_{n\hat{\zeta}_n}^{\tau}(\hat{\psi}_n)=0$, (\ref{l9.e1})~implies
$\sqrt{n}U_{0\zeta_0}^{\tau}(\hat{\psi}_n)=$
\[\sqrt{n}U_{0\hat{\zeta}_n}^{\tau}(\hat{\psi}_n)+o_P^{{\cal H}_1}(1)
=-\sqrt{n}\left(U_{n\zeta_0}^{\tau}(\psi_0)-U_{0\zeta_0}^{\tau}(\psi_0)\right)
+o_P^{{\cal H}_1}(1)=O_P^{{\cal H}_1}(1),\]
where $O_P^{B}(1)$ denotes a term bounded in probability
uniformly over the set $B$. By lemma~\ref{l4}, we know that there
exists a constant $e_2>0$ such that 
\[\|U_{0\zeta_0}^{\tau}(\psi)-
U_{0\zeta_0}^{\tau}(\psi_0)\|_{{\cal H}_1}\geq 
e_2\|\psi-\psi_0\|_{\infty}
+o(\|\psi-\psi_0\|_{\infty}),\]
as $\|\psi-\psi_0\|_{\infty}\rightarrow 0$.
Hence $\sqrt{n}\|\hat{\psi}_n-\psi_0\|_{\infty}(e_2-o_P(1))\leq O_P(1)$,
and we obtain the second conclusion of the lemma.

For the third part, we have 
\[\sqrt{n}\sup_{t\in[0,\tau]}\left|
\pp_n W(t;\hat{\theta}_n)-PW(t;\hat{\theta}_n)
\right|=\sqrt{n}\sup_{t\in[0,\tau]}|(\pp_n-P)W(t;\theta_0)|+o_P(1)\] 
$=O_P(1)$ and $\sqrt{n}\sup_{t\in[0,\tau]}|PW(t;\hat{\theta}_n)-
PW(t;\theta_0)|=O_P(1)$ by the first two parts of this lemma. 
Hence $\sqrt{n}\sup_{t\in[0,\tau]}
\left|\pp_n W(t;\hat{\theta}_n)-PW(t;\theta_0)\right|=O_P(1)$.
The result now follows by the Lipschitz continuity of $\log(u)$
over strictly positive compact intervals.$\Box$

{\it Proof of lemma~\ref{l10}.} 
The first inequality follows from the definitions.
For the second inequality, we use a Taylor's expansion around
$(\hat{\zeta}_n,\hat{\gamma}_n,\hat{\Gamma}_n)$ to obtain
$\tilde{X}_n(\hat{\zeta}_n,\hat{\lambda}_n)-\tilde{X}_n(\hat{\zeta}_n,
\lambda_0)=$
\[-\pp_n U_{\hat{\zeta}_n}^{\tau}(\hat{\gamma}_n,
A_n^{(\hat{\Gamma}_n)})(\lambda_0-\hat{\lambda}_n)
-\psi_{n,t}^{(\lambda_0-\hat{\lambda}_n)}\left(\pp_n
\hat{\sigma}_{
\left(\hat{\zeta}_n,\hat{\gamma}_{n,t},A_n^{(\hat{\Gamma}_{n,t})}\right)}
\right)(\lambda_0-\hat{\lambda}_n),\]
for some $t\in[0,1]$, where $\hat{\lambda}_{n,t}\equiv
(\hat{\gamma}_{n,t},\hat{\Gamma}_{n,t})$; $\hat{\gamma}_{n,t}\equiv
t\hat{\gamma}_n+(1-t)\gamma_0$; $\hat{\Gamma}_{n,t}\equiv
t\hat{\Gamma}_n+(1-t)\Gamma_0$; and,
for any $h\in{\cal H}_{\infty}$, 
$\psi_{n,t}^{(h)}\equiv\left(h_1,h_2,h_3,\int_0^{(\cdot)}h_4(s)
dA_n^{(\hat{\Gamma}_{n,t})}(s)\right)$.
The score term is zero by definition of the NPMLE, and the second term
has absolute value bounded by $\hat{K}_n 
\|\hat{\lambda}_n-\lambda_0\|_{\infty}^2$,
where $\hat{K}_n$ is bounded in probability 
by the uniform consistency of $\hat{\lambda}_n$ and by
the form of the information terms listed in section~5.2. 

Now, letting $\psi_n(\gamma,\Gamma)\equiv(\gamma,A_n^{(\Gamma)})$, we have
$\tilde{X}_n(\hat{\zeta}_n,\lambda_0)-\tilde{X}_n^{\ast}(\hat{\zeta}_n)$
\begin{eqnarray}
&&\label{l10.e1}\\
&&=\pp_n\left\{\left(\ind\{Y\leq\hat{\zeta}_n\}-\ind\{Y\leq\zeta_0\}\right)
\right.\nonumber\\
&&\left.\mbox{\hspace{0.1in}}
\times\left[l_1^{\psi_n(\gamma_0,\Gamma_0)}(V,\delta,Z)
-l_2^{\psi_n(\gamma_0,\Gamma_0)}(V,\delta,Z)
-l_1^{\psi_0}(V,\delta,Z)+l_2^{\psi_0}(V,\delta,Z)\right]\right\}\nonumber
\end{eqnarray}
$=\int_0^{\tau}\pp_n\left\{
\left(\ind\{Y\leq\hat{\zeta}_n\}-\ind\{Y\leq\zeta_0\}\right)
\tilde{Y}(s)\tilde{K}_n(s)\right\}e^{-\Gamma_0(s)}
\left[d\tilde{G}_n(s)-d\tilde{G}_0(s)\right]$,
where
\begin{eqnarray*}
\tilde{K}_n(s)&=&\left[\dot{G}(H_1^{\psi_{n,t}}(V))
-\delta\frac{\ddot{G}(H_1^{\psi_{n,t}}(V))}{\dot{G}(H_1^{\psi_{n,t}}(V))}
\right]e^{\beta_0'Z(s)}\\
&&-\left[\dot{G}(H_2^{\psi_{n,t}}(V))
-\delta\frac{\ddot{G}(H_2^{\psi_{n,t}}(V))}{\dot{G}(H_2^{\psi_{n,t}}(V))}
\right]e^{\beta_0'Z(s)+\alpha_0+\eta_0'Z_2(s)}
\end{eqnarray*}
and $\psi_{n,t}\equiv\left(\gamma,\int_0^{(\cdot)}\Gamma_0(u)\left[
td\tilde{G}_n(u)+(1-t)d\tilde{G}_0(u)\right]\right)$, for
some $t\in[0,1]$, by the mean value theorem.
By the conditions given in section~2, we have that there is a
constant $k^{\ast}<\infty$ such that 
$\|\tilde{K}_n(s)\Gamma_0(s)\|_{v}\leq k^{\ast}$
with probability~1 for all $n\geq 1$. Thus the absolute value
of~(\ref{l10.e1}) is bounded above by
$k^{\ast}\|\tilde{G}_n-\tilde{G}_0\|_{\infty}\times\pp_n
\left|\ind\{Y\leq\hat{\zeta}_n\}-\ind\{Y\leq\zeta_0\}\right|
=O_P(n^{-1})$.
This last statement follows because $\|\tilde{G}_n-\tilde{G}_0\|_{\infty}
=O_P(n^{-1/2})$, 
$(\pp_n-P)\left|\ind\{Y\leq\hat{\zeta}_n\}-\ind\{Y\leq\zeta_0\}\right|
=o_P(n^{-1/2})$, and $P\left|
\ind\{Y\leq\hat{\zeta}_n\}-\ind\{Y\leq\zeta_0\}\right|=O_P(n^{-1/2})$
by theorem~\ref{t.l9}. Now the desired result follows.$\Box$

{\it Proof of lemma~\ref{l11}.} Note first that
\[\tilde{D}_n(\zeta)=\sqrt{n}(\mathbb{P}_n-P)\left\{\left[
\ind\{Y\leq \zeta\}-\ind\{Y\leq\zeta_0\}\right]\times
\left[l_1^{\psi_0}-l_2^{\psi_0}\right](V,\delta,Z)\right\}.\]
Denote $\tilde{H}\equiv[l_1^{\psi_0}-l_2^{\psi_0}](V,\delta,Z)$,
and note that $|\tilde{H}|\leq c_{\ast}$ almost surely
for a fixed constant $c_{\ast}<\infty$. Thus 
$F_{\epsilon}\equiv\ind\{\zeta_0-\epsilon\leq Y
\leq\zeta_0+\epsilon\}c_{\ast}$ 
serves as an envelope for
the class of functions 
\[{\cal F}_{\epsilon}\equiv\{\left[\ind\{Y\leq\zeta\}
-\ind\{Y\leq\zeta_0\}\right]\tilde{H}:|\zeta-\zeta_0|\leq\epsilon\},\]
for each $\epsilon>0$.
Note that by the assumptions on the density $\tilde{h}$ in a neighborhood
of $\zeta_0$, we have for some $\epsilon_3>0$ that there exists
$0<k_{\ast},k_{\ast\ast}<\infty$ such that $k_{\ast}\epsilon\leq
\tilde{p}(\epsilon)\equiv
P[\zeta_0-\epsilon\leq Y\leq\zeta_0+\epsilon]\leq k_{\ast\ast}\epsilon$ 
for all $0\leq\epsilon\leq\epsilon_3$. 
Thus the bracketing entropy 
\[N_{[]}(u\|F_{\epsilon}\|_{P,2},{\cal F}_{\epsilon},L_2(P))\leq
O\left(\frac{\epsilon}{u^2\tilde{p}(\epsilon)}\right)\leq
O\left(\frac{1}{c_{\ast}u^2}\right),\]
for all $u>0$ and $0\leq\epsilon\leq\epsilon_3$; 
and thus, by theorem~2.14.2 of \cite{vw96}, 
there exists a $c_{\ast\ast}<\infty$
such that 
\[E\left[\sup_{|\zeta-\zeta_0|\leq\epsilon}
|\tilde{D}(\zeta)|\right]\leq c_{\ast\ast}\|F_{\epsilon}\|_{P,2}
\leq c_{\ast\ast}c_{\ast}\sqrt{k_{\ast\ast}\epsilon},\]
for all $0\leq\epsilon\leq\epsilon_3$.
The result now follows for $k_2=c_{\ast\ast}c_{\ast}
\sqrt{k_{\ast\ast}}$.$\Box$

{\it Proof of theorem~\ref{t3}.} We can deduce from section~3 that 
\begin{eqnarray*}
\lefteqn{\tilde{L}_n(\hat{\psi}_n,\zeta_{n,u})
-\tilde{L}_n(\hat{\psi}_n,\zeta_0)}
&&\\ 
&\mbox{\hspace{0.5cm}}=&
\mathbb{P}_n\left\{\left(\ind\{\zeta_{n,u}<Y\leq\zeta_0\}
-\ind\{\zeta_0<Y\leq\zeta_{n,u}\}\right)
\left[l_2^{\hat{\psi}_n}-l_1^{\hat{\psi}_n}\right](V,\delta,Z)\right\}\\
&\mbox{\hspace{0.5cm}}=&n^{-1}Q_n(u)+\hat{E}_n(u),\;\;\;\;\mbox{where}
\end{eqnarray*}
$\hat{E}_n(u)\equiv\mathbb{P}_n\left\{\left(\ind\{Y\leq\zeta_0\}
-\ind\{Y\leq\zeta_{n,u}\}\right)\left[l_2^{\hat{\psi}_n}-l_2^{\psi_0}
-l_1^{\hat{\psi}_n}+l_1^{\psi_0}\right](V,\delta,Z)\right\}$.
By arguments similar to those used in the proof of lemma~\ref{l10},
we can obtain constants $0<F_1,F_2<\infty$ such that
$\left|l_j^{\hat{\psi}_n}(V,\delta,Z)-l_j^{\psi_0}(V,\delta,Z)\right|
\leq F_j\|\hat{\psi}_n-\psi_0\|_{\infty}$ almost surely, for $j=1,2$. Hence
\[|\hat{E}_n(u)|\leq\pp_n\left|\ind\{Y\leq\zeta_0\}
-\ind\{Y\leq\zeta_{n,u}\}\right|O_P(n^{-1/2}).\]
By arguments given in the proof of lemma~\ref{l11}, we know that
\[(\pp_n-P)\left|\ind\{Y\leq\zeta_0\}
-\ind\{Y\leq\zeta_{n,u}\}\right|=O_P^{\mathbb{U}_{n,M}}(n^{-1}).\]
Since also $\sup_{u\in\mathbb{U}_{n,M}}P\left|\ind\{Y\leq\zeta_0\}
-\ind\{Y\leq\zeta_{n,u}\}\right|=O(n^{-1})$ by condition B2(i),
we now have that $\hat{E}_n=O_P^{\mathbb{U}_{n,M}}(n^{-3/2})$.
The desired result now follows.$\Box$

{\it Proof of theorem~\ref{t4}.} Fix $h\in{\cal H}_{\infty}$.
We first establish that
$\left(Q_n^+,{\cal Z}^n(h)\equiv\right.$ 
$\left.\sqrt{n}\pp_n U_{\zeta_0}^{\tau}(\psi_0)(h)\right)$
converges weakly to $(Q^+,{\cal Z}(h))$, on $D_M\times\re$,
where $Q^+$ and ${\cal Z}(h)$ are independent, for each 
fixed $M<\infty$, and ${\cal Z}(h)$ is mean zero Gaussian with
variance $\tilde{\sigma}_h^2\equiv\mbox{var}[U_{\zeta_0}^{\tau}(\psi_0)(h)]$. 
Accordingly, fix $M$, and 
let $0=u_0<u_1<u_2<\cdots<u_J\leq M$ be a finite collection of
points and $q_1,\ldots,q_J,\tilde{q}$ be arbitrary real numbers. Our plan
is to first show that the characteristic function of
$(Q_n^+(u_1),\ldots,Q_n^+(u_J),{\cal Z}^n(h))$ converges to that of
$(Q^+(u_1),\ldots,Q^+(u_J))$ times that of ${\cal Z}(h)$. 
Since the choice of points
$u_1,\ldots,u_J$ is arbitrary, this will imply convergence
of all finite-dimensional distributions. We will then show
that $Q_n^+$ is asymptotically tight, and this will imply
the desired weak convergence. 

Let $y\mapsto I_{nj}(y)\equiv\ind\{\zeta_0+u_{j-1}/n<y\leq\zeta_0+u_j/n\}$, 
$j=1,\ldots,J$; and 
$F_i\equiv[l_1^{\psi_0}-l_2^{\psi}](V_i,\delta_i,Z_i)$ and
${\cal Z}_i\equiv U_{\zeta_0}^{\tau}(\psi_0)(h)(X_i)$, 
$i=1,\ldots,n$. In other words, ${\cal Z}_i$ is the score contribution
from the $i$th observation. Thus
\begin{eqnarray}
\lefteqn{P\exp\left[i\left\{\sum_{j=1}^Jq_j[Q_n^+(u_j)-Q_n^+(u_{j-1})]
+\tilde{q}{\cal Z}^n(h)\right\}\right]}\mbox{\hspace{1.0in}}&&\label{t4.e1}\\
&=&\prod_{k=1}^n P\left[\exp\left\{\sum_{j=1}^J i q_jI_{nj}(Y_k)F_k\right\}
e^{i\tilde{q}{\cal Z}_k/\sqrt{n}}\right].\nonumber
\end{eqnarray}
However, using the facts that 
$e^{\sum_j w_j}-1=\sum_j(e^{w_j}-1)$ when 
only one of the $w_j$'s differs from zero and 
$e^{uv}-1=u(e^{v}-1)$ when $u$ is dichotomous, we have
$\exp\left\{\sum_{j=1}^J iq_jI_{nj}(Y_k)F_k\right\}
=1+\sum_{j=1}^J\left(e^{iq_j I_{nj}(Y_k)F_k}-1\right)
=1+\sum_{j=1}^JI_{nj}(Y_k)\left(e^{iq_jF_k}-1\right)$.
Combining this with condition~B2 and the boundedness of
$F_k$ and ${\cal Z}_k$, we obtain
$P\left[\exp\left\{\sum_{j=1}^J i q_jI_{nj}(Y_k)F_k\right\}
e^{i\tilde{q}{\cal Z}_k/\sqrt{n}}\right]$
\begin{eqnarray*}
&=&Pe^{i\tilde{q}{\cal Z}_k/\sqrt{n}}+
\sum_{j=1}^J\frac{(u_j-u_{j-1})\tilde{h}(\zeta_0)}{n}
P\left[\left.\left(e^{iq_jF_k}-1\right)e^{i\tilde{q}{\cal Z}_k/\sqrt{n}}
\right|Y=\zeta_0+\right]\\
&&+o(n^{-1})\\
&=&1+n^{-1}\left[-\frac{\tilde{q}^2\tilde{\sigma}_h^2}{2}
+\tilde{h}(\zeta_0)\sum_{j=1}^J(u_j-u_{j-1})\{\phi^+(q_j)-1\}\right]
+o(n^{-1}),
\end{eqnarray*}
where $o(1)$ denotes a quantity going to zero uniformly
over $k=1,\ldots,n$. Thus the right-hand side of~(\ref{t4.e1}) is
\[\exp\left[\frac{-\tilde{q}^2\tilde{\sigma}_h^2}{2}+
\tilde{h}(\zeta_0)\sum_{j=1}^J(u_j-u_{j-1})\{\phi^+(q_j)-1\}\right],\]
which is precisely 
\[P\exp\left[i\tilde{q}{\cal Z}(h)+i\sum_{j=1}^jq_j\left\{
Q^+(u_j)-Q^+(u_{j_1})\right\}\right].\]
Thus the finite dimensional distributions converge as desired.

We next need to verify that $Q_n^+$ is asymptotically tight
on~$[0,M]$. Since there exists
a constant $c_{\ast}<\infty$ such that $\max_{1\leq i\leq n}
|F_i|\leq c_{\ast}<\infty$ almost surely, we have that
$|Q_n^+(u_2)-Q_n^+(u_1)|\leq c_{\ast}n\pp_n
\ind\{\zeta_0+u_1/n<Y\leq\zeta_0+u_2/n\}$,
for all $0\leq u_1<u_2\leq M$. Thus we are done if we can show
that $u\mapsto\tilde{R}_n(u)\equiv n\pp_n\ind\{\zeta_0<Y\leq\zeta_0+u/n\}$ is
tight on $[0,M]$. To this end, fix $0\leq u_1<u_2\leq M$. Now,
the expectation of $|\tilde{R}_n(u_2)-\tilde{R}_n(u_1)|$ is
$nP\{\zeta_0+u_1/n<Y\leq\zeta_0+u_2/n\}\rightarrow 
|u_2-u_1|\tilde{h}(\zeta_0)$, as $n\rightarrow\infty$. 
This implies the desired tightness since 
$u\mapsto\tilde{R}_n(u)$ is monotone. We have now established that
$\left(Q_n^+,{\cal Z}^n(h)\right)$
converges weakly to $(Q^+,{\cal Z}(h))$, on $D_M\times\re$,
where $Q^+$ and ${\cal Z}(h)$ are independent, for each 
fixed $M<\infty$. Similar arguments also yield the weak convergence
of $\left(Q_n^-,{\cal Z}^n(h)\right)$ to
$(Q^-,{\cal Z}(h))$, on $D_M\times\re$,
where $Q^-$ and ${\cal Z}(h)$ are again independent, for each 
fixed $M<\infty$. Thus also $\left(Q_n,{\cal Z}^n(h)\right)$ converges
weakly to $(Q,{\cal Z}(h))$, on $D_M\times\re$,
where $Q$ and ${\cal Z}(h)$ are independent, for each 
fixed $M<\infty$. Since $n(\hat{\zeta}_n-\zeta_0)=O_P(1)$, the
argmax continuous mapping theorem (theorem~3.2.2 of \cite{vw96}) now yields
that $\left(n(\hat{\zeta}_n-\zeta_0),{\cal Z}^n(h)\right)$
converges weakly to $\left(\argmax\,Q,{\cal Z}(h)\right)$, with
the desired asymptotic independence. The remaining
results follow.$\Box$

{\it Proof of theorem~\ref{t5}.} We have
\begin{eqnarray*}
0&=&\sqrt{n}\pp_n U_{\hat{\zeta}_n}^{\tau}(\hat{\psi}_n)\\
&=&\sqrt{n}\pp_n U_{\zeta_0}^{\tau}(\hat{\psi}_n)+
\sqrt{n}(\pp_n-P)\left(U_{\hat{\zeta}_n}^{\tau}(\hat{\psi}_n)
-U_{\zeta_0}^{\tau}(\hat{\psi}_n)\right)\\
&&+\sqrt{n}
P\left(U_{\hat{\zeta}_n}^{\tau}(\hat{\psi}_n)
-U_{\zeta_0}^{\tau}(\hat{\psi}_n)\right)\\
&\equiv&\sqrt{n}\pp_n U_{\zeta_0}^{\tau}(\hat{\psi}_n)+B_{1,n}+B_{2,n},
\end{eqnarray*}
where the index set for the score terms is ${\cal H}_1$.
By arguments similar to those used in the proof of theorem~\ref{t.l9},
combined with the fact that $n(\hat{\zeta}_n-\zeta_0)=O_P(1)$, we have
that both $B_{1,n}=o_P^{{\cal H}_1}(1)$ and $B_{2,n}=o_P^{{\cal H}_1}(1)$. 
Thus $\sqrt{n}\pp_n U_{\zeta_0}(\hat{\psi}_n)=o_P^{{\cal H}_1}(1)$. 
We also have that
\[\sqrt{n}(\pp_n-P)U_{\zeta_0}^{\tau}(\hat{\psi}_n)-
\sqrt{n}(\pp_n-P)U_{\zeta_0}^{\tau}(\psi_0)=o_P^{{\cal H}_1}(1).\]
Combining this with lemma~\ref{l4}, the Z-estimator master theorem
(theorem~3.3.1 of \cite{vw96}) now yields the desired results.$\Box$

{\it Proof of corollary~\ref{c1}.} We first derive the unconditional
limiting distribution of $\sqrt{n}(\hat{\psi}_n^{\circ}-\psi_0)$.
If a class of measurable functions
${\cal F}$ is $P$-Glivenko-Cantelli with $\|P\|_{\cal F}<\infty$, then
the class $\kappa\cdot{\cal F}=\{\kappa f:f\in{\cal F}\}$, where
$\kappa$ denotes a generic version of one of the weights $\kappa_i$,
is also $P$-Glivenko-Cantelli, by theorem~3 of \cite{vw00}.  
Thus we can apply the
results of theorem~\ref{t1}, with only minor modification, combined
with the simple fact that $\bar{\kappa}\rightarrow\mu_{\kappa}$
almost surely, to yield that $\hat{\psi}_n^{\circ}\rightarrow\psi_0$
outer almost surely. Note that the proof is made somewhat easier than
before since we already know $\hat{\zeta}_n\rightarrow\zeta_0$
almost surely. Furthermore, if a class of measurable functions~${\cal F}$
is $P$-Donsker with $\|P\|_{\cal F}<\infty$, then the multiplier central
limit theorem (theorem~2.9.2 of \cite{vw96})
yields that the class $\kappa\cdot{\cal F}$ is also $P$-Donsker.
Hence we can apply the results of theorem~\ref{t4}, with only minor
modification, to yield that $\sqrt{n}(\hat{\psi}_n^{\circ}-\psi_0)$
is asymptotically linear with influence function
$\tilde{l}^{\circ}(h)=(\kappa/\mu_{\kappa})U_{\zeta_0}^{\tau}
(\sigma_{\theta_0}^{-1}(h))$, $h\in{\cal H}_1$. The factor 
$\mu_{\kappa}^{-1}$ occurs because the information operator for
the weighted version of the likelihood is $\mu_{\kappa}\sigma_{\theta_0}$.
We now have that $\sqrt{n}(\hat{\psi}_n^{\circ}-\hat{\psi}_n)
=\sqrt{n}\pp_n(\kappa/\mu_{\kappa}-1)U_{\zeta_0}^{\tau}
(\sigma_{\theta_0}^{-1}(\cdot))+o_P^{{\cal H}_1}(1)$, unconditionally. 

Finally,  
the conditional multiplier central limit theorem (theorem~2.9.6 
of \cite{vw96})
yields part~(ii) of the theorem. The factor $(\mu_{\kappa}/\sigma_{\kappa})$
arises because $\mbox{var}(\kappa/\mu_{\kappa})=\sigma_{\kappa}^2
/\mu_{\kappa}^2$. Similar arguments establish~(i) by using
parallel Glivenko-Cantelli and Donsker results for the nonparametric
bootstrapped empirical process.$\Box$

{\it Proof of lemma~\ref{s9.l1}.} Let $\mu(x)$ denote the baseline 
measure and $\rho_n(x)$, $\rho(x)$ the density function under 
$P_n$ and $P$ respectively. In the general situation, 
verifying~(\ref{c8.e1}) is equivalent to finding a function $h$
such that:
\begin{eqnarray*}
\;\;\;\;\lefteqn{\int 
\left[ \frac{\left(\frac{dP_n(x)}{d \mu(x)}\right)^{1/2 }-
\left(\frac{dP(x)}{d \mu(x)}\right)^{1/2}}{1/\sqrt{n}} -
\frac{1}{2}h(x)\left(\frac{dP(x)}{d \mu(x)}\right)^{1/2}\right]^2d\mu(x)}&&\\
&=& \int \left[ \frac{\rho_n(x)^{1/2}-\rho(x)^{1/2} }{1/\sqrt{n}}
- \frac{1}{2}h(x)\rho(x)^{1/2} \right]^2 d \mu(x)\\
& \rightarrow & \int \left [ \frac{1}{2} 
\frac{\dot {\rho}(x)}{(\rho(x))^{1/2}}-
\frac{1}{2}h(x)\frac{\rho(x)}{(\rho(x))^{1/2}} \right]^2 d\mu(x)\\
&=&\int \left [ \frac{1}{2} \frac{\dot{\rho}(x)}
{\rho(x)}(\rho(x))^{1/2}
-\frac{1}{2}h(x)(\rho(x))^{1/2} \right]^2 d\mu(x)\\
&=&0.
\end{eqnarray*}
Hence the given score function satisfies~(\ref{c8.e1})
by the smoothness of the log-likelihood.$\Box$

{\it Proof of lemma~\ref{s9.l2}.} Note that a consequence of the 
Donsker theorem for contiguous alternatives
(theorem~3.10.12 of \cite{vw96}) 
is that for any bounded $P$-Donsker class ${\cal F}$,
$\|\pp_n-P\|_{\cal F}\weakpn 0$. Thus the proof of
lemma~\ref{l2} can be reconstituted to yield 
that $\|\hat{A}_0\|_{[0,\tau]}$ is bounded in probability
under $P_n$, since all of the classes of functions involved are bounded
$P$-Donsker classes. We can similarly modify the proof of theorem~\ref{t1}
to yield the desired results since, once again, the only classes of
functions involved are bounded and $P$-Donsker. This is true,
in particular, for the key class given in lemma~\ref{l.t1.1}, for
any $k<\infty$. Thus $\|\hat{\psi}_0-\psi_0^{\ast}\|_{\infty}
\weakpn 0$.$\Box$

{\it Proof of theorem~\ref{s9.t1}.} The basic idea of the proof
is to use the Donsker theorem for contiguous alternatives in
combination with key arguments in the proof of theorem~\ref{t5}
and the form of the score and information operators under
model C2'. Pursuing this course, we obtain for any $(h_1,h_2)\in
\re^{q+1}$,
\begin{eqnarray*}
(h_1,h_2')\hat{S}_1(\zeta)&=&\sqrt{n}\pp_n(1,1)\left[\left(
\begin{array}{c}U_{\zeta,1}^{\tau}\\ U_{\zeta,2}^{\tau}
\end{array}\right)(\psi_0^{\ast})\left(\begin{array}{c}
h_1\\h_2\end{array}\right)\right.\\
&&\left.-\left(
\begin{array}{c}U_{\zeta_0,3}^{\tau}\\ U_{\zeta_0,4}^{\tau}
\end{array}\right)(\psi_0^{\ast})\left([\sigma_{\ast}^{22}]^{-1}
\sigma_{\ast}^{21}(\zeta)\left(\begin{array}{c}h_1\\h_2
\end{array}\right)\right)\right]+o_{P_n}^{[a,b]}(1)\\
&\equiv&\sqrt{n}\pp_n H_{\ast}(\zeta)+o_{P_n}^{[a,b]}(1),
\end{eqnarray*}
where $o_{P_n}^B(1)$ denotes a quantity going to zero in
probability, under $P_n$, uniformly over the set $B$. Now
the Donsker theorem for contiguous alternatives yields that
the right-hand side converges to a tight, Gaussian process with
covariance $P[H_{\ast}(\zeta_1)H_{\ast}(\zeta_2)]$, for all
$\zeta_1,\zeta_2\in[a,b]$, and mean $P\left[H_{\ast}\left\{
U_{\zeta_0,1}^{\tau}(\psi_0^{\ast})(\alpha_{\ast})+U_{\zeta_0,2}^{\tau}
(\psi_0^{\ast})(\eta_{\ast})\right\}\right]$. Note that we
only need to compute the moments under the null distribution~$P$.
Careful calculations verify that this yields the desired results.$\Box$

{\it Proof of corollary~\ref{c2}.} The limiting results under
$P_n$ follow from theorem~\ref{s9.t1} and the continuous mapping
theorem, provided we can show that 
\begin{eqnarray}
\inf_{\zeta\in[a,b],v\in\re^{q+1}:\|v\|=1}v'V_{\ast}(\zeta)v&>&0.
\label{e1.m9}
\end{eqnarray}
The limiting null distribution 
results will similarly follow from the
fact that under the null distribution~$P$, $\nu_{\ast}(\zeta)=0$
for all $\zeta\in[a,b]$. Note that in both the null and alternative
settings, $V_{\ast}(\zeta)$ only depends on the null limiting
distribution. It is sufficient to verify that $\sigma_{\psi_0^{\ast},
\zeta_n}$ is one-to-one for all sequences $\zeta_n\in[a,b]$
and $h_n\in{\cal H}_{\infty}$. Note that we can ignore any
differences between $\zeta_0$ and $\zeta$ in calculating
$\zeta\mapsto\sigma_{\psi_0^{\ast},\zeta}^{22}$ because of
the non-identifiability of $\zeta$ under the null
hypothesis, ie., $\zeta\mapsto\sigma_{\psi_0^{\ast},\zeta}^{22}$
is constant. Assume now that there exists sequences $\zeta_n\in[a,b]$
and $h_n\in{\cal H}_{\infty}$ such that 
$\sigma_{\psi_0^{\ast},\zeta_n}h_n\rightarrow 0$. We will now
show that this forces $h_n\rightarrow 0$. Without loss of
generality, we can assume $\zeta_n\rightarrow\zeta_{\ast}$ and
$h_n\rightarrow h$. Since the map $h\mapsto\sigma_{\psi_0^{\ast},\zeta}h$
is continuous and since $\zeta\mapsto\sigma_{\psi_0^{\ast},\zeta}h$ is
cadlag, we can further assume without loss of generality that
either $\sigma_{\psi_0^{\ast},\zeta_{\ast}}h=0$ or that
$\sigma_{\psi_0^{\ast},\zeta_{\ast}^-}h=0$ (the $\zeta_{\ast}^-$ denotes
that we are converging to $\zeta_{\ast}$ from below). The arguments for either
case are the same, so we will for brevity only give the proof for
the first case.

By the arguments surrounding expressions~(\ref{c12:e9}), (\ref{c12:e11}) 
and~(\ref{c12:e12}), combined with the non-identifiability of $\zeta$
under the null model, we obtain that expression~(\ref{c12:e12}) must now
hold for all $t\in(0,\tau]$ but with $\zeta_{\ast}$ replacing $\zeta_0$.
In ortherwords,
$\tilde{Y}(t)(h_1\ind(Y>\zeta_{\ast})+
h_2'Z_2(t)\ind(Y>\zeta_{\ast})+h_3'Z+h_4(t))=0$, almost surely,
for all $t\in(0,\tau]$.
Since var$[Z(t_4)|Y>\zeta_{\ast}]\geq\mbox{var}[Z(t_4)|Y>b]
\times\pr{Y>b}/\pr{Y>\zeta_{\ast}}$ 
is positive definite by condition~B4, we have $h_3=0$.
We can similarly use~B4 to verify that var$[Z(t_3)|Y\leq\zeta_{\ast}]$
is positive definite and thus $h_2=0$. Now $h_1=0$ and $h_4=0$
easily follow. Hence $h\mapsto\sigma_{\psi_0^{\ast},\zeta}h$ is
uniformly one-to-one in a manner 
which yields the conclusion~(\ref{e1.m9}).$\Box$

{\it Proof of theorem~\ref{s9.t2}.} The results follow from 
arguments similar to those used in the proof of theorem~\ref{s9.t1},
but based on the conditional multiplier central limit theorem
for contiguous alternatives, theorem~\ref{s9.t2.t1} below.$\Box$

\begin{theorem}\label{s9.t2.t1} (Conditional multiplier central
limit theorem for contiguous alternatives) Let ${\cal F}$ be
a $P$-Donsker class of measurable functions, and let $P_n$ satisfy
\[\int\left[\sqrt{n}(dP_n^{1/2}-dP^{1/2})-\frac{1}{2}hdP^{1/2}
\right]^{1/2}\rightarrow 0\label{s9.t2.t1.e1},\]
as $n\rightarrow\infty$, for some real valued, measurable
function $h$. Also assume
$\lim_{M\rightarrow\infty}$ $\limsup_{n\rightarrow\infty}
P_n(f-Pf)^2\ind\{|f-Pf|>M\}=0$
for all $f\in{\cal F}$, and that the multipliers in
the weighted bootstrap, $\kappa_1,\ldots,\kappa_n$, are i.i.d.
and independent of the data, with mean $0<\mu_{\kappa}<\infty$
and variance $0<\sigma_{\kappa}^2<\infty$, and with
$\int_0^{\infty}\sqrt{P(\kappa_1>u)}du<\infty$. 
Then $(\mu_{\kappa}/\sigma_{\kappa})(\pp_n^{\circ}-\pp_n)
\weakpnboot\mathbb{G}$ in $\ell^{\infty}({\cal F})$,
where $\mathbb{G}$ is a tight, mean zero Brownian bridge process.
\end{theorem}

{\it Proof.} The detailed proof can be found in chapter~11
of Kosorok (To appear). We now present a synopsis of the proof.
Let $\tilde{\kappa}_i\equiv\sigma_{\kappa}^{-1}
(\kappa_i-\mu_{\kappa})$, $i=1,\ldots,n$, and note that
\begin{eqnarray}
\label{s9.t2.t1.e3}&&\\
\pp_n^{\circ}-\pp_n&=&n^{-1/2}\sum_{i=1}^n(\kappa_i/\bar{\kappa}
-1)\Delta_{X_i}\;\;=\;\;
n^{-1/2}\sum_{i=1}^n(\kappa_i/\bar{\kappa}-1)(\Delta_{X_i}-P)
\nonumber\\
&=&\frac{\sigma_{\kappa}}{\mu_{\kappa}}n^{-1/2}\sum_{i=1}^n
\tilde{\kappa}_i(\Delta_{X_i}-P)+
\left(\frac{\sigma_{\kappa}}{\bar{\kappa}}-\frac{\sigma_{\kappa}}
{\mu_{\kappa}}\right)n^{-1/2}\sum_{i=1}^n\tilde{\kappa}_i(\Delta_{X_i}-P)
\nonumber\\
&&+\left(\frac{\mu_{\kappa}}{\bar{\kappa}}-1\right)n^{-1/2}
\sum_{i=1}^n(\Delta_{X_i}-P),\nonumber
\end{eqnarray}
where $\Delta_{X_i}$ is the Dirac measure of the observation $X_i$.
Since ${\cal F}$ is $P$-Donsker, we also have that
$\dot{\cal F}\equiv\{f-Pf:f\in{\cal F}\}$ is $P$-Donsker. Thus
by the unconditional multiplier central limit theorem, 
we have that $\tilde{\kappa}\cdot{\cal F}$ is also $P$-Donsker. Now, by
that fact that $\|P(f-Pf)\|_{\cal F}=0$ (trivially) combined with
the central limit theorem under contiguous alternatives, we have that
both $f\mapsto n^{-1/2}\sum_{i=1}^n\tilde{\kappa}_i(\Delta_{X_i}-P)f
\weakpn\mathbb{G}f$ and $f\mapsto n^{-1/2}\sum_{i=1}^n(\Delta_{X_i}
-P)\weakpn\mathbb{G}f+P[(f-Pf)h]$ in $\ell^{\infty}({\cal F})$.
Thus the last two terms in~(\ref{s9.t2.t1.e3})$\;\weakpn 0$, and hence
$\sqrt{n}(\mu_{\kappa}/\sigma_{\kappa})(\pp_n^{\circ}-\pp_n)
\weakpn\mathbb{G}$ in $\ell^{\infty}({\cal F})$. 
This now implies the unconditional asymptotic
tightness and desired asymptotic measurability of $\sqrt{n}
(\mu_{\kappa}/\sigma_{\kappa})(\pp_n^{\circ}-\pp_n)$.
Fairly standard arguments can now be used along with the given pointwise
uniform square integrability condition to verify that 
$\sqrt{n}(\mu_{\kappa}/\sigma_{\kappa})(\pp_n^{\circ}-\pp_n)$
applied to any finite dimensional collection $f_1,\ldots,f_m\in{\cal F}$
converges under $P_n$ in distribution, conditional on the data,
to the appropriate limiting Gaussian process. This now implies
$\sqrt{n}(\mu_{\kappa}/\sigma_{\kappa})(\pp_n^{\circ}-\pp_n)\,
\weakpnboot\mathbb{G}$.$\Box$

{\it Proof of corollary~\ref{c3}.} Assume at first
that $\tilde{M}_n$ is a fixed
number $\tilde{M}<\infty$. Theorem~\ref{s9.t2} now yields that the
collection $\{\hat{S}_{1,1}^{\circ}-\hat{S}_1,\ldots,
\hat{S}_{1,\tilde{M}_n}^{\circ}-\hat{S}_1\}$ converges jointly, conditionally
on the data, to $\tilde{M}$ i.i.d. copies of $\mathbb{Z}_{\ast}$.
Thus $\hat{V}_n$ converges weakly to the sample covariance
process (divided by $\tilde{M}_n$ instead of $\tilde{M}_n-1$)
of an i.i.d. sample of $\tilde{M}_n$ copies of
$\mathbb{Z}_{\ast}$. The same result holds true if we allow
$\tilde{M}_n$ to go to~$\infty$ slowly enough. Since the
Gaussian processes involved are tight, $\hat{V}_n$ will thus
be consistent for $\Sigma_{\ast}$, uniformly over $\zeta\in[a,b]$.
Similar arguments yield pointwise consistency of $\hat{\mathbb{F}}$
and $\tilde{\mathbb{F}}$ at continuity points of 
$\hat{\mathbb{T}}_{\ast}$ and $\tilde{\mathbb{T}}_{\ast}$. 
Since it is not hard to verify that
both $\hat{\mathbb{T}}_{\ast}$ and $\tilde{\mathbb{T}}_{\ast}$
have continuous distributions, the pointwise consistency extends
to the desired uniform consistency.$\Box$

\section*{Acknowledgments}
The authors thank Editor Morris Eaton, an associate editor, and two
referees for their extremely
careful review and helpful suggestions that led
to an improved paper.


\begin{thebibliography}{9}

\bibitem{a01}
{\sc Andrews, D. W. K.} (2001). Testing when a parameter is on the
boundary of the maintained hypothesis. {\it Econometrica} {\bf
69}, 683--73.

\bibitem{ap94}
{\sc Andrews, D. W. K., and Plogerger, W.} (1994). Optimal
tests when a nuisance parameter is present only under the
alternative. {\em Econometrica} {\bf 62}, 1383--1414.

\bibitem{bn04}{\sc Bagdonavi\v{c}ius, V., and Nikulin, M.} (2004).
Statistical modeling in survival analysis and its influence on the
duration analysis. {\it Advances in survival analysis}, 411--429,
{\it Handbook of Statistics, 23}. Elsevier, Amsterdam.

\bibitem{bkrw98}
{\sc Bickel, P. J., Klaassen, C. A. J., Ritov, Y. and Wellner,
J. A.} (1998). {\it Efficient and Adaptive Estimation for
Semiparametric Models}. Springer-Verlag, New York.

\bibitem{bd81}{\sc Bickel, P. J., and Doksum, K. A.} (1981). An analysis of
transformations revisited. {\em Journal of the American
Statistical Association} {\bf 76}, 296--311.

\bibitem{br97}{\sc Bickel, P. J., and Ritov, Y.} (1997). Local asymptotic
normality of ranks and covariates in transformation models.
{\em Festschrift for Lucien Le Cam: Research papers in probability
and statistics}, 43--54.

\bibitem{bc64}{\sc Box, G. E. P., and Cox, D. R.} (1964). 
An analysis of transformations. (With discussion) 
{\em Journal of the Royal Statistical Society}, Series B 
{\bf 26}, 211--252.

\bibitem{bc82}{\sc Box, G. E. P., and Cox, D. R.} (1982). 
An analysis of transformations revisited, rebutted. 
{\em Journal of the American Statistical Association}
{\bf 77}, 209--210.

\bibitem{c89}{\sc Chappell, R.} (1989). Fitting bent lines to data, with 
applications to allometry. {\it Journal of Theoretical Biology} 
{\bf 138}, 235-256.

\bibitem{cwy95}{\sc Cheng, S. C., Wei, L. J., and Ying, Z.} (1995). Analysis
of transformation models with censored data. {\em Biometrika} {\bf 82},
835--845.

\bibitem{cwy97}{\sc Cheng, S. C., Wei, L. J., and Ying, Z.} (1997). Predicting
survival probabilities with semiparametric transformation models. {\it Journal
of the American Statistical Association} {\bf 92}, 227--235.

\bibitem{dd88}{\sc Dabrowska, D.M. and Doksum, K.A.} (1988). Estimation
and Testing in the Two-sample Generalized Odds-Rate Model. {\it
Journal of the American Statistical Association} {\bf 83}, 1--23.

\bibitem{d87}{\sc Davies, R. B.} (1987). Hypothesis testing when a nuisance
parameter is present only under the alternative. {\em Biometrika}
{\bf 74}, 33--43.

\bibitem{fyw98}
{\sc Fine, J. P., Ying, Z., and Wei, L. J.} (1998). On the linear
transformation model for censored data. {\em Biometrika} {\bf 85}, 980--986.

\bibitem{ih81}{\sc Ibragimov, I. A., and Has'minskii, R. Z.} (1981). {\em
Statistical estimation: Asymptotical theory}. Springer, New York.

\bibitem{kta}{\sc Kosorok, M. R.} (To appear). {\em Introduction to Empirical
Processes and Semiparametric Inference}. Springer, New York.

\bibitem{klf04}{\sc Kosorok, M. R., Lee, B. L. and Fine, J. P.} (2004). Robust
Inference for Univariate Proportional Hazards Frailty Regression
Models. {\it The Annals of Statistics} {\bf 32}, 1448-1491.

\bibitem{lsl90}{\sc Liang, K.-Y., Self, S. G., and Liu, X.} (1990). The Cox
proportional hazards model with change point: An epidemiologic
application. {\em Biometrics} {\bf 46}, 783--793.

\bibitem{ly93}{\sc Lin, D. Y. and Ying, Z.} (1993). Cox regression with
incomplete covariate measurements. {\it Journal of the American Statistical
Association} {\bf 88}, 1341--1349.

\bibitem{lb97}{\sc Luo, X. and Boyett, J. M.} (1997). Estimation of a
threshold parameter in cox regression. {\it Communication in
Statistics--Theory and Methods} {\bf 26}, 2329--2346.

\bibitem{ltc97}{\sc Luo, X., Turnbull, B.W. and Clark, L.C.} (1997).
Likelihood ratio tests for a changepoint with survival data. {\it
Biometrica} {\bf 84}, 555--565.

\bibitem{mrv97} 
{\sc Murphy, S. A., Rossini, A. J., and van der Vaart, A. W.} (1997). 
Maximum likelihood estimation in the
proportional odds model. {\it Journal of the American Statistical
Association} {\bf 92}, 968--976.

\bibitem{p98}{\sc Parner, E.} (1998). Asymptotic theory for the correlated 
gamma-frailty model. {\em Annals of Statistics} {\bf 26}, 183--214.

\bibitem{p82}{\sc Pettit, A. N.} (1982). Inference for the linear model using
a likelihood based on ranks. {\em Journal of the Royal Statistical
Society}, Series B {\bf 44}, 234--243.

\bibitem{p84}{\sc Pettit, A. N.} (1984). Proportional odds models for survival
data and estimates using ranks. {\em Applied Statistics} {\bf 33}, 169--175.

\bibitem{pr94}{\sc Politis, D. N., and Romano, J. P.} (1994). Large sample
confidence regions based on subsamples under minimal assumptions.
{\em Annals of Statistics} {\bf 22}, 2031--2050.

\bibitem{p03}{\sc Pons, O.} (2003). Estimation in a cox regression model
with a change-point according to a threshold in a covariate. {\it
The Annals of Statistics} {\bf 31}, 442--463.

\bibitem{stg98}
{\sc Scharfstein, D. O., Tsiatis, A. A., and Gilbert, P. B.} (1998).
Semiparametric efficient estimation in the generalized odds-rate class
of regression models for right-censored time-to-event data.
{\em Lifetime Data Analysis} {\bf 4}, 355--391.

\bibitem{s98}{\sc Shen, X.} (1998). Proportional odds regression and sieve
maximum likelihood estimation. {\em Biometrika} {\bf 85}, 165--177.

\bibitem{sv04}{\sc Slud, E. V., and Vonta, F.} (2004). Consistency of the
NPML estimator in the right-censored transformation model.
{\em Scandinavian Journal of Statistics} {\bf 31}, 21--41.

\bibitem{v98}{\sc van der Vaart, A. W.} (1998). {\em Asymptotic Statistics}.
Cambridge University Press, Cambridge.

\bibitem{vw96}{\sc van der Vaart, A. W., and Wellner, J. A.} (1996). {\it
Weak Convergence and Empirical Processes: With Applications to
Statistics.} Springer, New York.

\bibitem{vw00}
{\sc van der Vaart, A. W., and Wellner, J. A.} (2000). Preservation
theorems for Glivenko-Cantelli and Uniform Glivenko-Cantelli classes.
{\em High Dimensional Probability II}, 113--132. Birkhauser, Boston.

\end{thebibliography}
\end{document}